\documentstyle[12pt,a4,amscd,amssymb,verbatim]{amsart}

\newtheorem{thm}{Theorem}

\newtheorem{lemma}{Lemma}

\begin{document}

\large
\title{\bf Grading the translation functors in type A } 
\author{Steen Ryom-Hansen} 
\address{Matematisk Afdeling, Universitetsparken 5 \\
DK-2100 K{\o}benhavn \O, Danmark }
%\address{University of Freiburg \\
%Eckerstra\ss e 1, 79104 Freiburg }
\begin{abstract}
{We deal with the representation theory of quantum 
groups and Hecke algebras at roots of unity.
We relate the philosophy of Andersen, Jantzen and Soergel 
on graded translated functors to the Lascoux, Leclerc and 
Thibon-algorithm. This goes via the Murphy standard 
basis theory and the idempotents coming from the
Murphy-Jucys operators. Our results lead to a guess on a 
tilting algorithm outside 
the lowest $p^2$ alcove, which at least in the 
$SL_2$-case coincides with Erdmann's results.}

\end{abstract}
\maketitle

\section {\bf Introduction }
\noindent

This work is concerned with the modular representation theory of the 
symmetric group and of the Hecke algebra. The main sources of inspiration 
are the work of Lascoux, Leclerc 
and Thibon [LLT] on the crystal basis of the Fock module and the paper 
of Andersen, Jantzen and Soergel [AJS]
on the representation theory of Frobenius 
kernels.  

\medskip 
\noindent

The philosophy of [AJS] is to provide the representation theory --
including the Jantzen translation functors -- with 
a grading, that explains the fact that the Kazhdan-Lusztig algorithm constructs
polynomials rather than numbers.
Under the Schur functor, (certain of) the translation functors
correspond to the classical $ r $--inducing and
$ r $--restricting operators from the modular representation theory of the 
symmetric group. On the other hand, the 
action of the quantum group $ U_q(\widehat{{\mathfrak sl}_l}) $
on the Fock space can be viewed as $q$-deformations of these 
operators on the Grothendieck group level.
So according to the above 
philosophy, it is natural to look for a representation theoretical meaning
of these deformed operators.

\medskip 
\noindent

We here propose to connect them to the invariant form on
the Specht modules. This is quite a different approach than that 
of [AJS]. Our calculations are considerably simpler than those of [AJS], on 
the other hand we have to pass to the characteristic zero situation on 
the way and therefore we are not able to construct a 
grading on the representation category itself. 

\medskip 
\noindent

We still believe that our results and methods 
provides new insight to modular representation 
theory. For example, we show 
how our results lead to a guess on a tilting algorithm outside 
the lowest $p^2$ alcove, which at least in the 
$SL_2$-case coincides with Erdmann's results.

\medskip 
\noindent

The paper is organized as follows. In the next
section we review briefly the representation theory of 
the symmetric group, following James's book [J]. 
In the third section we perform the calculations 
for the restriction functors. In the fourth section 
we treat the induction functors. These are the 
hardest calculations of the paper and, as a matter 
of fact, we were forced to go via the 
representation theory of the Hecke algebra 
where we can rely on Murphy's standard basis.
In the fifth section we treat the duality.
In the last section we present the $SL_2$-algorithm.

\section{\bf Preliminaries}

We will use the terminology of James book [J]. 
So let $ k $ be an arbitrary field and 
let $ S_n $ denote the
symmetric group on $ n $ letters acting on $ \{1,\ldots,n \} $ on the right.
Let $ \lambda = ( \lambda_1, \lambda_2, \ldots , \lambda_r ) $ be 
a partition of $ n $ into $ r $ parts and let $ S_{\lambda} $ be the 
corresponding Young 
subgroup of $ S_n $ which is defined as the row stabilizer of the 
$\lambda $-tableau 
$ t^{\lambda} $ in which the numbers $ 1,2,\ldots,n $ are placed in the
diagram along the rows. 
The permutation module $ M^{\lambda} $ 
is obtained by inducing the trivial representation $ k $ from $ S_{\lambda} $ 
to $ S_n $. It has a basis consisting of so called $ \lambda $-tabloids; 
these are equivalence classes of $ \lambda $-tableaux, where 
two tableaux $ t_1 $ and $ t_2 $ are said to be equivalent if there is an 
element $ \pi $ of the row stabilizer $ R_{t_1} $ of $ t_1 $, such that 
$ t_1  \pi = t_2 $. We denote the tabloid class of the tableau $ t $ by 
$ \{t\} $.

\medskip
(Although we shall usually assume that $ \lambda $ is 
a partition, the above 
constructions and statements would also make sense for compositions, 
i.e. for unordered partitions.
For instance, if $ V $ is an $ r $-dimensional vector space over $k$, 
we can make 
$  V^{ \otimes n  } $ into an $ S_n $-module by place permutation and as 
such we then have that 
$$ V^{ \otimes n } = \bigoplus_{\lambda   } 
M^{\lambda } $$ 
with $ \lambda $ running over all compositions $ \lambda $ of $ n$ 
in less than $ r $ parts. 
In other words, there is a tabloid basis 
of $ V^{ \otimes r} $ as well.) 

\medskip

Let $ C_t $ be the column stabilizer of the 
$\lambda$-tableau $ t $. We associate then to $ t $ the element $ e_t \in 
M(\lambda) $ as follows: 
$$ e_t = \sum_{\pi \in C_t }{\mbox{sign}} \pi  \{ t \pi \} $$
The Specht module $ S( \lambda ) $ labeled by $ \lambda $ is defined as 
the subspace of $ M( \lambda ) $ generated by the $ e_t $'s with $ t $ running 
over the set of all $ \lambda $-tableaux. 

\section{\bf Restriction Functors }

In this section we take a careful look at the proof of the branching rule 
given in James's book [J]. According to this the restricted module 
$ {\rm{Res}}^{S_n}_{S_{n-1}}\, S(\lambda) $ has a Specht filtration, in which 
$ S(\mu ) $ 
occurs as a subquotient if and only if $ \mu $ is a subdiagram of 
$ \lambda $ such that 
the difference $ \lambda \setminus \mu $ consists of exactly one node.

\medskip
\noindent

This filtration can be constructed in terms of the standard basis 
$ \{ \,e_t \, |\, \, \, t \mbox{ standard tableau}\,  \} $ 
of $ S( \lambda ) $ 
in the following way: 
Let $ r_1 <  r_2 <
\ldots < r_m $ be the row numbers of the rows 
from which a node can be removed 
to leave a diagram. Denote by $ V_i $ the subspace of $ S(\lambda ) $ 
generated by the $ e_t $'s where $n$ is in the $ r_1 $'th, $ r_2 $'th, 
$ \ldots $ or $ r_i $'th row of $ t $. Then we have  
a filtration of $S_{n-1}$-modules
$$ V_1 \subseteq
V_2 \subseteq \ldots \subseteq  V_m, \, \, \, \, \,  V_i / V_{i-1} 
\cong S(\lambda^i ) $$
with $ \lambda^i $ denoting the partition that is obtained from 
$ \lambda $ by
removing the last node from the $ r_i$'th row. 
This construction is independent of the characteristic of $k$.

\medskip
\noindent

Recall the $ S_n $-invariant form 
$ \langle \cdot,\cdot \rangle _{\lambda} $ 
on the 
permutation module $ M(\lambda) $ which is defined by the formula
$$ \langle \{ t_1\},\{ t_2 \} \rangle_{\lambda} \,= 
\, \delta_{\{\, t_1\},
\{ t_2\}} $$
Its restriction to the Specht module $S(\lambda) \subset M(\lambda) $ 
is also denoted $ \langle \cdot,\cdot \rangle _{\lambda}$.

\medskip

The following theme is central to the paper: 
Assume that $ k = \mathbb Q $. Then 
$ \langle \cdot,\cdot \rangle_{\lambda} $ is nondegenerate.
Let $ U_i $ be the 
complement with respect to $ \langle \cdot,\cdot \rangle _{\lambda} $ of 
$ V_{i-1} $ in $ V_i $. We then have that $$ U_i \stackrel{\iota}{\cong}
V_i \, /\, V_{i-1}
\stackrel{\pi}{\cong} S(\lambda^i) $$
Here $ \iota   $ is induced from the 
injection of $ U_i \subseteq V_i $, while 
$ \pi $ is the map that takes $ e_t \in S(\lambda) $ to 
$ e_{t^{\prime}} \in S(\lambda^i) $, where $t^{\prime}$ is obtained 
from $ t $ by removing the node containing $ n $.
We can now define a new scalar product 
$ \langle \cdot,\cdot \rangle _{\lambda}^{\prime} $
on $ U_i $ by pulling back the standard product 
$ \langle \cdot,\cdot \rangle _{\lambda^i}$ on $ S( \lambda^i ) $ to 
$ U_i $. Since $ U_i \, \cong \,  S(\lambda^i) $ is simple when 
$ \mbox{char}k=0 $, we have that  
$ \langle \cdot,\cdot \rangle_{\lambda} $ and 
$ \langle \cdot,\cdot \rangle_{\lambda}^{\prime} $ only differ by a 
scalar in $k$ which we denote by $m_i^2$: 
$$ \langle \cdot,\cdot \rangle_{\lambda} = m_i^2 
\langle \cdot,\cdot \rangle_{\lambda}^{\prime} $$
\noindent
James and Murphy [JM] calculated the scalars $ m_i^2 $.
The expression provides the inductive step in their determination of the Gram 
matrix and so they call it the branching 
rule for the Gram matrix. It goes as follows: 
\begin{equation}
m_i^2 \,  = \, \prod^{r_i-1}_{a=1} \, \frac{H_{ac}}{H_{ac}-1} 
\end{equation}
where $ ( r_i,c ) $ are the coordinates of the 
last node of the $r_i$'th row we are looking at and 
$ H_{ac}$
denotes the hook length of the hook centered at the node $ ( a,c ) $.

\medskip 
\noindent
{\bf Example:} 
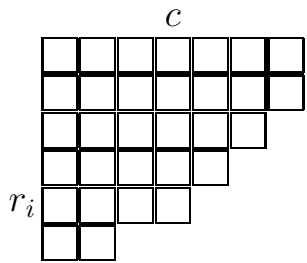
\unitlength0.5cm
\begin{center}
\begin{picture}(7,8)

\multiput(0,1)(1,0){2}{\framebox(0.9,0.9){}}
\multiput(0,2)(1,0){3}{\framebox(0.9,0.9){}}
\multiput(3,2)(1,0){1}{\framebox(0.9,0.9){}}
\multiput(0,3)(1,0){3}{\framebox(0.9,0.9){}}
\multiput(3,3)(1,0){1}{\framebox(0.9,0.9){}}
\multiput(4,3)(1,0){1}{\framebox(0.9,0.9){}}
\multiput(0,4)(1,0){3}{\framebox(0.9,0.9){}}
\multiput(3,4)(1,0){3}{\framebox(0.9,0.9){}}
\multiput(0,5)(1,0){7}{\framebox(0.9,0.9){}}
\multiput(0,6)(1,0){7}{\framebox(0.9,0.9){}}
\put(3,7){\makebox(0.9,0.9){$c$}}

\put(-1,2){\makebox(0.9,0.9){$r_i$}}

\end{picture}
\end{center}
Here we find $ m_i^2 \,= \, \frac{8}{7} \,
\frac{7}{6}\, \frac{5}{4} \, \frac{3}{2} \,= \, \frac{5}{2} $.

\medskip
\noindent
Let us fix $ l \in \mathbb N $.
In the LLT setup which we consider next, $l$ is the order 
of the root of unity.
We are especially interested in 
the ``quantum $l$-adic valuation'' of the numbers $m^2_i$ 
which we define as follows 
%$$ \nu_{q}(m^2_i) = \left\{ \begin{array}{ll}
%                           1 & \mbox{ if $ l \, | \,m $ }  \\
%                           0 & \mbox{ otherwise } 
%                          \end{array}
%                 \right. $$

\begin{equation}\label{valuation} \nu_{q}(m^2_i) = \nu_{\frak p}
\left( \prod^{r_i-1}_{a=1} \, \frac{[H_{ac}]_v}{[H_{ac}-1]_v} \right)
\end{equation}
where $ [n]_v = \frac{v^{n+1}-v^{-n-1}}{v-v^{-1}}$ is the usual 
Gaussian integer, while $\nu_{\frak p}: {\mathbb Z}[v,v^{-1}] \mapsto
\mathbb Z $ is the valuation with 
respect to the $l$'th cyclotomic polynomial.
Notice that if all the hook lengths $H_{ij} $ are less than 
$l$, this ``quantum valuation'' coincides with the
usual $l$-adic valuation of the number $m_i^2$.

\noindent

\medskip

Let us consider the $ q $-analogue of the Fock space 
$ {\cal F}_q $ as 
introduced by Hayashi. It is defined as 
$$ {\cal F}_q = \bigoplus_\lambda 
{\mathbb Q}( q ) \, | \lambda \rangle $$ with basis parameterized 
by the set of 
all partitions $ \lambda $. There is an action of the quantum group 
$ U_q(\widehat{{\mathfrak sl}_l}) $ on $ {\cal F}_q $ ($ q $ is here 
transcendent), making $ {\cal F}_q $ into an integrable module. Using this, 
Lascoux, Leclerc and Thibon defined a kind of
Kazhdan-Lusztig algorithm, which first of all calculates the 
global crystal basis of $ {\cal F}_q $, but also, by a theorem of Ariki,
the decomposition numbers of the representation theory of the Hecke algebra
at an $l$'th root of unity.
We are aiming at a connection between this algorithm and the numbers 
$ \nu_{{\frak p}}(m_i^2 ) $
coming from the branching rule. Let us therefore give the precise 
formulas 
for the action of $ U_q(\widehat{{\mathfrak sl}_l}) $ on $ {\cal F}_q $.

\medskip

In order to do this we need the concept of an $i$-node of the Young 
diagram $ Y(\lambda) $ of $\lambda$.
We start by filling in the nodes of $Y(\lambda)$ with integers, 
increasingly 
along the rows, decreasingly along the columns and starting with $ 0 $ in 
the $ (1,1) $'th position. We then reduce these numbers mod $ l $, the 
resulting diagram is called the $l$-diagram of the partition. We call a 
node an $ i $-node $ i \in {\mathbb Z} / l \, \mathbb Z $ if is has been 
filled in with $ i $.
A node of the diagram is called removable if it can be removed
to leave another diagram, and a (virtual) node that can be added to 
the 
diagram to yield another Young diagram, is called an indent node. 
Finally, we 
define for a subdiagram $ \mu $ of $ \lambda $ such that 
$ \gamma \, = \, 
\lambda \setminus \mu $ consist of just one node,
$I_i^r(\lambda, \mu )$ (resp. $ R_i^r(\lambda, \mu ) $) as the 
number of indent (removable) $ i $-nodes situated to the right of 
$ \gamma $.
Let now $ f_i, k_i $ and $ e_i $ be the standard 
generators of the quantum group $ U_q(\widehat{\mathfrak{sl}_l}) $. 
Define 
\begin{equation} N_i^r(\lambda, \mu ) = I_i^r(\lambda, \mu )  - 
R_i^r(\lambda, \mu ),\,\,\,\,\,\,
N_i^l(\lambda, \mu )  = I_i^l(\lambda, \mu )  - 
R_i^l(\lambda, \mu ) 
\end{equation}
Then $ f_i $ and $ e_i $ operate on  
$ {\cal F}_q $ by the following formulas
\begin{equation}\label{action}  
f_i  \lambda = \sum_{\mu} \, 
q^{N_i^r(\lambda, \mu ) }\, | \mu \rangle \,\,\,\,\,\,
e_i  \mu = \sum_{\lambda} \, 
q^{-N_i^l(\lambda, \mu ) }\, | \lambda \rangle
\end{equation}
where the first sum is over all partitions
$ \mu $ such that $ \lambda \setminus \mu $ consists of one $i$-node 
and analogously in the second expression.
There are similar formulas for the 
action of the other generators $ k_i $ and $ D $ of $U_q(\widehat{{\frak
 sl}_l}) $, see eg. [LLT].

We can now formulate the following Theorem.
\begin{thm} The numbers $ N_i^r(\lambda, \mu ) $ and 
$ \nu_q(  m_i^2  ) $ coincide.
\end{thm}
\begin{pf*}{Proof} 
Let us consider the $ (i,j ) $ hook of the $ l $-diagram. 
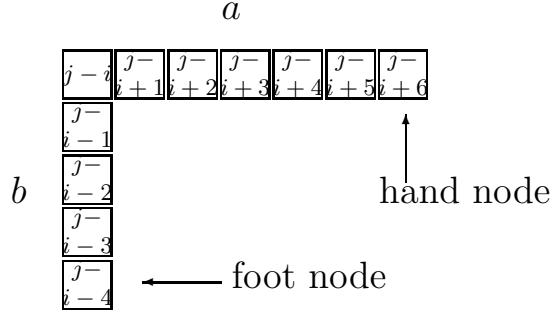
\unitlength0.7cm
\begin{center}
\begin{picture}(7,7)

\makebox{
\scriptsize
\multiput(0,1)(1,0){1}{\framebox(0.9,0.9){$\begin{array}{c} j- \\ i-4
\end{array}$}}
\multiput(0,2)(1,0){1}{\framebox(0.9,0.9){$\begin{array}{c} j- \\ i-3
\end{array}$}}
\multiput(0,3)(1,0){1}{\framebox(0.9,0.9){$\begin{array}{c} j- \\ i-2
\end{array}$}}
\multiput(0,4)(1,0){1}{\framebox(0.9,0.9){$\begin{array}{c} j- \\ i-1
\end{array}$}}
\multiput(0,5)(1,0){1}{\framebox(0.9,0.9){$\begin{array}{c} j- i
\end{array}$}}
\multiput(1,5)(1,0){1}{\framebox(0.9,0.9){$\begin{array}{c} j-  \\ i+1
\end{array}$}}
\multiput(2,5)(1,0){1}{\framebox(0.9,0.9){$\begin{array}{c} j- \\ i+2
\end{array}$}}
\multiput(3,5)(1,0){1}{\framebox(0.9,0.9){$\begin{array}{c} j- \\ i+3
\end{array}$}}
\multiput(4,5)(1,0){1}{\framebox(0.9,0.9){$\begin{array}{c} j- \\ i+4
\end{array}$}}
\multiput(5,5)(1,0){1}{\framebox(0.9,0.9){$\begin{array}{c} j- \\ i+5
\end{array}$}}
\multiput(6,5)(1,0){1}{\framebox(0.9,0.9){$\begin{array}{c} j- \\ i+6
\end{array}$}}}
\multiput(3,6.5)(1,0){1}{$a$}
\multiput(-1,3)(1,0){1}{$b$}
\put(3,1.5){\vector(-1,0){1.5}}
\put(3.2,1.4){foot node}
\put(6.5,3.4){\vector(0,1){1.3}}
\put(6,3){hand node}
\end{picture}
\end{center}

Here $a$ denotes the arm length of the hook, while $ b$ denotes the foot
length of the hook. We now have that $$ \begin{array}{c}H_{ij} = a+b-1 
\equiv (\, hand -(j-i)
+1\, ) + (\, j-i - foot + 1 \,)  \\ \equiv hand -foot + 1 
\mod l  \end{array} $$
\noindent
with $ hand $ (resp. $ foot $) denoting the residue of the hand node 
(foot node) of the hook. 

\medskip
Let us now calculate the constant $ \nu_q( m_i^2) $ at the 
partition
$ \lambda $. We can wlog. assume that the $ (\,i,l\,) $'th node is 
a $ 0 $-node.
The $ l$-diagram will then have the following form:

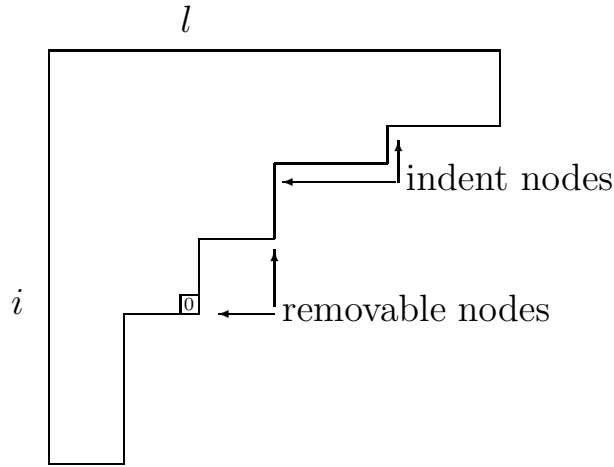
\unitlength0.5cm
\begin{center}
\begin{picture}(13,13)
\makebox{
\line(0,1){11}
\put(0,11){\line(1,0){12}}
\put(12,11){\line(0,-1){2}}
\put(12,9){\line(-1,0){3}}
\put(9,9){\line(0,-1){1}}
\put(9,8){\line(-1,0){3}}
\put(6,8){\line(0,-1){2}}
\put(6,6){\line(-1,0){2}}
\put(4,6){\line(0,-1){2}}
\put(4,4){\line(-1,0){2}}
\put(2,4){\line(0,-1){4}}
\put(2,0){\line(-1,0){2}}
\put(3.5,4){\line(0,1){0.5}}
\put(3.5,4.5){\line(1,0){0.5}}
\put(-1,4){$ i $ }
\put(3.5,11.5){$ l $}
\makebox{ 
\tiny 
\put(3.6,4.1){$ 0 $} }
\put(6,4){\vector(-1,0){1.5}}
\put(6.2,3.8){removable nodes}
\put(6,4.2){\vector(0,1){1.5}}
\put(9.2,7.5){\vector(-1,0){3}}
\put(9.5,7.3){indent nodes}
\put(9.3,7.5){\vector(0,1){1.1}}
}
\end{picture}

\end{center}

Let $ \lambda^{\prime} > \lambda_i $ be some row length of the diagram and 
let $ s, t $ be such that $$ \{ s, s+1 , \ldots , t  \} = 
\{  j  \mid \lambda_j = \lambda^{\prime}  \} $$
Then the contribution to $ \nu_q(  m_i^2  ) $ coming from 
$ \{ s, s+1 , \ldots , t  \} $ equals 
$$ \nu_{\frak p} \left( \frac{[H_{tl}]_v}{[H_{sl}-1]_v} \right)$$

Now, according to the above calculation $ l $ divides the hook 
length $ H_{il}$
if and only if the last node of the $ i $´th row is a $ l-1 $ node. On 
the other hand, this means that there is an indent $ 0 $ node in the
position to the right of 
this row. Analogously, we have that $ l  \mid  H_{il} -1 $ if and only 
if the last node of the $ i $´th row is a $ 0 $-node: thus a removable 
$ 0 $-node. We conclude that the $t$'th row will contribute (with $ 1 $) 
to  $ \nu_q(  m_i^2  ) $ if and only if it gives rise to an indent node, 
while the 
$ t $'th row will contribute (with $ -1 $) to $ \nu_q( m_i^2  ) $, if 
and only its end node is a removable one. 

\medskip
It just remains to consider the rows of length $ \lambda_i $. So let 
$ i^{\prime} $ be given by 
$$ \{ i, i+1 , \ldots , i^{\prime}  \}  = 
\{  j  \mid \lambda_j = \lambda_i  \} $$
But then the corresponding rows will contribute with $ H_{i^{\prime}l } $ to
$ m_i^2 $. And once again $ l $ divides the hook length if and only if 
there is an indent $ 0 $-node in the position beyond the $ i^{\prime} $'th row, so 
also here things match up. We have proved the Theorem.  
\end{pf*}
{\bf Remark:}
The action of $e_i$ is a deformation of Robinson's 
$r$-restriction functors, so we would have expected 
$ \nu_q(  m_i^2  ) =
-N_i^l(\lambda, \mu ) $ 
in the above Theorem to get the complete match with 
(\ref{action}). Let us therefore consider the action $e^l_i$ (resp. $f^l_i$) on
${\cal F}_q $ which is equal to that of $e_i$ (resp. $f_i$) but with 
$ q^{- N_i^l(\lambda, \mu )} $ 
(resp. $q^{ N_i^r(\lambda, \mu )} $) replaced by $q^{ N_i^r(\lambda,\mu )} $
(resp. $q^{- N_i^l(\lambda, \mu )} $).
Define furthermore 
$$ \varphi: {\cal F}_q \rightarrow  {\cal F}_q
\,\,\,\,\,\,\, \lambda \mapsto  q^{[\lambda]_l} \lambda $$
where $ [\lambda]_l $ is the $ l$-weight of the partition
$ \lambda $, i.e. the number of $l$-skewhooks one should 
remove to arrive at the $l$-core. Then we have that 
$$ \varphi \circ e_i^l = e_i  \circ \varphi, \,\,\,\,\,\,
\varphi \circ f_i^l = f_i  \circ \varphi $$
which comes from the formula, [LLT, equation 13]
$$ -[\lambda]_l = -[\mu]_l + N_i^l(\lambda, \mu ) +
N_i^r(\lambda, \mu ) $$
It is in other words basically the same problem to 
determine the lower global basis $ G^l(\lambda) $ with respect to
$ f_i^l,e_i^l$ as with respect to $ f_i,e_i$.
The relation is 
\begin{equation}
G^l(\lambda) = q^{[\lambda]_l} G(\lambda)
\end{equation}

\section{ \bf Induction functors. }

We wish in this section to prove that the $ f_i $-operators on the Fock
space can be realized 
in a similar way as the $ e_i $-operators. This is of course a conclusion 
one might expect. Actually the LLT algorithm only uses the 
$ f_i $-operators so from our point of view, this 
conclusion is more important than the one 
for the $ e_i $-operators.

\medskip
\noindent
  
A first approach towards this result might be to use a 
Frobenius reciprocity argument. We were however unable to find 
any arguments along that line.
A reason why this approach apparently does not work 
may be that not all functors admit adjoints.

\medskip
\noindent

A major difference between the induction functors and the restriction 
functors
is that, unlike the restricted module, the induced module does not 
come with a
natural basis which is compatible with the Specht filtration, so to obtain 
results for the induction functors one
cannot just copy the calculations done for the restriction functors. 
Now, it is possible to extract a basis of the induced Specht module 
from the 
chapter on the Littlewood Richardson rule in James book [J], which indeed is 
compatible with the Specht filtration -- this basis is however not 
well suited for the diagonalization process -- the calculations
explode very quickly.
 
\medskip
\noindent

We choose to work with the 
Hecke algebra setting. Although this seems like a further 
complication of 
the problem, it provides us with a more natural setting for the 
simultaneous 
induction and diagonalization. The reason for this is first of all 
Murphy's standard basis which 
turns out to behave well with respect to the diagonalization process.
The results on the symmetric group 
case can then be obtained by specializing the Hecke algebra parameter 
$ q $ to $ 1 $.

\medskip
\noindent

Let us now review the basic definitions of the representation theory of the 
Hecke algebra $ H_n $ of type A; contrary to the symmetric group case the 
Specht modules are here constructed as ideals in $ H_n $ itself.

\medskip
\noindent

Let the standard generators of the Hecke algebra $ H_n $ of type A be 
$ T_w , \, w \in S_n $; they satisfy the relations 
\begin{equation} T_w \, T_v  = \left\{ \begin{array}{ll}
                           T_{wv} & \mbox{ if $ l(wv)\, = \, l(w)\,+\,1 $ }  \\
                           q \,  T_{wv} \, + \, ( q - 1 ) \, T_w & 
                           \mbox{ otherwise } 
                          \end{array}
                 \right. 
\end{equation}

\begin{equation}
 T_v \, T_w  = \left\{ \begin{array}{ll}
                           T_{vw} & \mbox{ if $ l(vw)\, = \, l(w)\,+\,1 $ }  \\
                           q \,  T_{vw} \, + \, ( q - 1 ) \, T_w & 
                           \mbox{ otherwise } 
                          \end{array}
                 \right. 
\end{equation}
\noindent
with $ l(\cdot ) $ denoting the standard length function on $ S_n $. One
checks that these relations make $ H_n $ into an associative algebra 
over 
$ {\mathbb Z}[ q, q^{-1} ] $ with unit element $ T_1 $ and with
$ \{ \, T_w
 \mid  w \in S_n \, \} $ as a basis.

\medskip 
\noindent
Now for $ X \, \subseteq \, S_n $ we define the following elements of 
$ H_n $
$$ \iota( X ) \, = \, \sum_{ w \in X }\, T_w, \, \, \, \mbox{and} \, \, \
\epsilon( X ) \, = \, \sum_{ w \in X }\,( -q)^{-l(w)} T_w $$
Let $h^*$ denote the image of $ h \, \in H_n $ under the 
antiautomorphism of $ H_n $ induced by the map $ T_w  \mapsto 
T_{w^{-1}},\, w 
\in S_n $. For any row standard $ \lambda $-tableau $ t $ we define 
$ d(t) 
\in S_n $ as the element of $ d  \in S_n $ satisfying 
$ t \,= \, t^{\lambda} \, d $. We denote the row stabilizer of 
$ t^{\lambda} $ by $ S_{\lambda} $ and can introduce the 
following elements of $ H_n $
$$ x_{st} \,= \, T^*_{d(s)} \, \iota( S_{\lambda} )\, T_{ d(t) } 
\mbox{ and  }
   y_{st} \,= \, T^*_{d(s)} \, \epsilon( S_{\lambda} )\, T_{ d(t) } $$
\noindent
With $ \lambda $ running through all partitions of $ n $ and $ s, t $ 
through all row standard 
$ \lambda $-tableaux, the $ x_{st} $ (as well as the $ y_{st} $) 
form a 
basis of the Hecke algebra: the Murphy standard basis of $ H_n $. 
It is a cellular basis in the sense of Graham and Lehrer [GL], indeed
it is in many senses the {\it prototype} of a cellular basis. 

\medskip
\noindent

Let $ s $ and $ t^{\prime} $ be row standard $ \lambda $ and 
$ \lambda^{\prime} $
tableaux, respectively; then we define $$ z_{st} \,= 
\, T^*_{d(s)}\, x_{\lambda
\lambda}\, T_{w_{\lambda}} \, y_{\lambda^{\prime} \lambda^{\prime}} 
\,  T_{ d(t^{\prime}) } $$
where the subscript $ \lambda $ stands for the tableau $ t^{\lambda} 
$ and 
where $ w_{\lambda} \in S_n $ is defined by the property that the 
elements of 
$ t^{\lambda}  w_{\lambda}  $ are entered by columns. The Specht module 
$ S(  \lambda  ) $ is now the right ideal of $ H_n $ generated by
$ z_{\lambda \lambda } $, it has a basis consisting of 
$$ \{ \,  z_{\lambda t }\, 
\mid t \mbox{  is a standard $ \lambda $ tableau } \}. $$
%For all of this one should consult Murphy´s paper [M1] or the more  
%axiomatic cellular basis approach [GL].
\medskip
\noindent
Our first task will be to describe a basis of 
$ \mbox{Ind}\,S(\lambda) $: the $H_n$-Specht module 
induced up to $H_{n+1}$. 
Let then $ \lambda^e $ be the partition of $ n  + 1 $ which 
is equal to $ \lambda $ except that $ ( \lambda^e )_1 = \lambda_1 +1 $. 
Dually we introduce the partition $ \lambda_e $ of $ n+1 $ by $\lambda_e
\, = \,( (\lambda^{\prime})^e)^{\prime} $. For 
a pair of partitions $(\lambda, \mu)$ we write 
$(\lambda, \mu)^{\prime} :=(\lambda^{\prime}, \mu^{\prime})$
and similarly for pairs of tableaux.
We were
unable to find the following Lemma in the literature:
\medskip
\begin{lemma}
$ \mbox{ \rm Ind}\,S(\lambda) \, = \, \mbox{span} \{\, 
x_{\lambda \lambda}\, T_{w_{\lambda}} \, y_{(\lambda^e, s)^{\prime}}
\, \mid \, s \, \,\,\, \lambda^e\mbox{-tableau } \} $
\end{lemma} 
\begin{pf*}{Proof} 
By the definition we have that $$ S(\lambda) = \mbox{span}\, \{\, 
x_{\lambda \lambda}\, T_{w_{\lambda}} \, y_{\lambda^{\prime} s^{\prime}}
\, \mid \, s \, \,\,\, \lambda \mbox{-tableau }\, \} $$
Now the transpositions $ \{ (1,n+1),(2,n+1),\ldots,(n,n+1) \} $ 
form a set 
of coset representatives of 
$ S_{n}\, \backslash \, S_{n+1} $ and we have
$$ \mbox {Ind}\, S(\lambda ) \, = 
\, S(\lambda) \otimes_{H_n} H_{n+1} \, = \, 
\bigoplus_{j=1}^n  \, S(\lambda)\,
T_{(j,n+1)}  $$ as $ {\mathbb Z}[ q, q^{-1} ] $-modules. But then the 
fact that 
$ T_{(j,n+1)} $ is invertible implies that $ \mbox {Ind}\, S(\lambda ) $
is free over $ {\mathbb Z}[ q, q^{-1} ] $. We can also work out the 
rank by passing to the quotient field $ {\mathbb Q}( q ) $ and using 
Frobenius reciprocity: we know everything about the restriction functors. The 
rank is given by the branching rule.

\medskip 

Now it is clear from the definition that $ y_{\lambda^{\prime}
\lambda^{\prime }} \,
= \, y_{(\lambda^{\prime})_e (\lambda^{\prime})_e}
= \, y_{(\lambda^e,\lambda^e)^{\prime}} $ 
and then
$$ x_{\lambda \lambda}\, T_{w_{\lambda}} \,
y_{\lambda^{\prime}\lambda^{\prime }} \, = \, 
x_{\lambda \lambda}\, T_{w_{\lambda}} \, 
y_{(\lambda^e,\lambda^e)^{\prime}} $$

Combining these equations we see that the right hand side of the Lemma
is a quotient of the left hand side. But by the next Lemma the 
$ H_{n+1} $-module $$\mbox{span} \{\, 
x_{\lambda \lambda}\, T_{w_{\lambda}} \, y_{(\lambda^e ,s)^{\prime}}
\, \mid \, s \, \,\,\, \lambda^e\mbox{-tableau }\, \}$$ has a Specht filtration
and is thus free over $ {\mathbb Z}[ q, q^{-1} ] $ as well. Since also 
this rank is given by the branching rule we are done.
\end{pf*}

Recall the total order $<$ on partitions introduced eg. in [M1]. See also
{\em loc. cit.} for the extension of this
ordering to tableaux, of arbitrary shapes, which we
shall also denote by $<$. For a partition $ \lambda $, the
tableau $ t^{\lambda} $ dominates all tableaux of shape $ \lambda $.
We then have the next Lemma:
\medskip

\begin{lemma} The $ H_{n+1}$-module $I = \mbox{span} \{\, 
x_{\lambda \lambda}\, T_{w_{\lambda}} \, y_{(\lambda^e, s)^{\prime}}
\, \mid \, s \, \,\,\, \lambda^e\mbox{-tableau } \} $ has a Specht filtration
$$ 0 \, \subseteq \, V_1 \, \subseteq V_2 \, \subseteq \, \ldots \,  V_m \, = 
\, I 
\,\, {\rm{with}} \,\,  V_i \, /\, 
V_{i-1} 
\cong S(\lambda^i ) $$ where the $ \lambda^i $´s are the partitions that
can be obtained from $ \lambda $ by adding exactly one node to $ \lambda $ 
ordered such that 
$ \lambda^1 \, \leq \, \lambda^2 \, \leq  \ldots \,\leq \lambda^m $.
\end{lemma} 

\medskip
\noindent

\begin{pf*}{Proof} Using Lemma 3.5 in [M1] (the Garnir relations) 
one can express
$y_{(\lambda^e, s)^{\prime}}  $ for all 
$\lambda^e$-tableaux
$ s $ as a linear 
combination of the elements
\begin{equation}\label{elements} \{ \, y_{u^{\prime}v^{\prime}} \, \mid \, u,v \mbox{ standard }\mu 
\mbox{-tableau such that } \mu \leq \lambda^e \} \end{equation}
On the other hand it is a general fact (see the remarks 
before Lemma 3.6 of [M1]) that for any $ w \, \in S_n $ and 
any partitions $ \tau , \nu $ of $ n$ we have
\begin{equation}\label{generalfact}
x_{\tau \tau }\,T_w\, y_{\nu \nu } 
\neq 0 \, \Rightarrow \, \tau \leq \nu^{\prime} 
\end{equation}
We have that $ x_{\lambda \lambda}\, = \, x_{\lambda_e \lambda_e} $ and
$  y_{u^{\prime} v^{\prime}} \, = \, T^*_{d(u^{\prime})}\,
 y_{\mu^{\prime} \mu^{\prime}}\,  T_{d(v^{\prime})} $ 
so we can use this general fact on the partitions $\lambda_e $ and $\mu^{\prime} $ 
of $n+1$. 
We find that 
only $ \mu $ with 
$ \lambda_e \, \leq \, \mu $,
can make 
$ x_{\lambda \lambda }\,T_{w_{\lambda}}\, y_{u^{\prime} v^{\prime}} $,
nonzero,
when $y_{u^{\prime} v^{\prime}}$ is from the set (\ref{elements}). 
But the only $\mu$ that satisfy the simultaneous 
inequalities $$ \lambda_e \, \leq \, \mu \, \leq \, \lambda^e $$
are those that can be 
obtained from $ \lambda $ by adding exactly one node.

\medskip 

We need a more precise analysis of the $ \mu $-tableaux
$ u $ that can occur as first index in 
$ y_{u^{\prime} v^{\prime}} $ from (\ref{elements}), in such a way that 
$ x_{\lambda \lambda }\,T_{w_{\lambda}}\, y_{u^{\prime} v^{\prime}} $
is nonzero.
These must first of all  
satisfy $ u^{\prime} \, \geq \, 
t^{(\lambda^e)^{\prime}} $; this is once again according to 
Lemma 3.5 in  [M1]. Let $ u^{\prime }\{\widehat{\,n+1\,}\} $ denote the 
tableau that is obtained from $ u^{\prime} $ by removing the $ n+1 $-node.
Then by definition of the order on tableaux we have 
$ u^{\prime }\{\widehat{\,n+1\,}\} \, \geq \,t^{\lambda^{\prime}}  $ or 
equivalently $ u\{\widehat{\,n+1\,}\} \, \leq \,t_{\lambda}  $. But then also
the underlying partition $ P \left( u\{\widehat{\,n+1\,}\} \right) $ of the 
tableau 
$  u\{\widehat{n+1}\} $ satisfies that
$$ P \left( u\{\widehat{n+1}\} \right) \, \leq \, \lambda. $$
On the other hand, since 
$ x_{\lambda \lambda}\, T_{w_{\lambda}} \, y_{u^{\prime} v^{\prime}} \, = 
 x_{\lambda_e \lambda_e} T_{w_{\lambda}} 
\, y_{u^{\prime} v^{\prime}} \
\neq \, 0 $ we have (as before) that $$ \lambda_e \, \leq \, P(  u )  $$ 
Since $ \mu $ is obtained from $ \lambda $ by adding one node, this is only 
possible if 
$$ P \left( u\{\widehat{n+1}\} \right) \, = \, \lambda $$
and since $ t_{\lambda} $ is the minimal tableau of shape $ \lambda $ we 
conclude from 
$ u\{\widehat{n+1}\} \, \leq \,t_{\lambda}  $
that in fact $$ u\{\widehat{n+1}\} \, = \, t_{\lambda}.  $$

\medskip
\noindent
So all in all we have proved that the only tableaux $ u $ that can occur 
as the 
first index of $ y_{u^{\prime}v^{\prime}} $ are those that are obtained from 
$ t_{\lambda} $ by adding an extra node which must contain the number $ n+1 $.
 
\medskip 
\noindent

%We must now show that all of these partitions are actually needed. To do so 
%recall the operators $ L_m, \, 1 \, \leq \, m \, \leq \, n $ on $ H_n $
%introduced by Murphy: 
%$$ L_m := q^{-1}\,T_{(m-1,m)} + q^{-2}\,T_{(m-2,m)} + \ldots +
% q^{-m}\,T_{(1,m)} $$

%Let $ \mu $ be a partition of $ n-1 $ with $ \overline{k} < k $ parts. Let 
%$ \mu_{i,e} $ be the tableau obtained from $ t^{\mu} $ by adjoining a node 
%containing $ n $ to the $ i $th row ( if possible to leave a tableau ). 
%Then Lemma 4.4 of [M1] yields that 
%$$  y_{\mu \mu} \left( - L_n  + [k]_q \right) = \sum_{i=1}^k \,
%q^{\alpha(i)} \,y_{\mu_{i,e}\mu_{i,e}} $$
%for some function $ \alpha $ ( we have applied the antiautomorphism $ * $ 
%to the actual formula of [ M1 ] ).

%Using this formula and induction we see that 
%at least the elements 
%$$ x_{\lambda \lambda}\, T_{w_{\lambda}} \, y_{\lambda_{i,e}^{\prime}
%\lambda_{i,e}^{\prime}} $$ will occur in the induced module. But then 
%for all $  \lambda_{i,e} $-tableaux $ s $ 
%$$ x_{\lambda \lambda}\, T_{w_{\lambda}} \, y_{\lambda_{i,e}^{\prime},
%s^{\prime}} $$ will also occur 
%by the operation with $ T_{d(s^{\prime})} $.

We must now show that all these partitions are actually needed. Let $ \mu $ 
be a partition of $ n-1 $ and let $ t $ be the tableau obtained from 
$ t^{\mu} $ by attaching $ n $ to the $ i $th row. Define 
$$ \alpha(\lambda)= \sum_i\,i\lambda_i $$ Then one has the formula (see 
the proof of Lemma 4.4 in [M1]).
\begin{equation}\label{Murphyformula}
 \begin{array}{l}
q^{\alpha(t)-l(d(t))}\,y_{tt}\, = \,
 q^{i} \left( 1- \sum_{c=j}^{n-1} \,
q^{c-m}\, T_{(c,n)} \right) q^{\alpha(\mu)}\,y_{\mu \mu}
+ \\ q^{i-1}\sum_{c=n}^{n-1} \,q^{c-n}\, T_{(c,n)}\, q^{\alpha(\mu)}\,
y_{\mu \mu} 
\end{array}
\end{equation}
where $ \{j,\ldots, m \} $ are the nodes of the $i$th row of $P(t)$.
Notice that there seems to be a, for our purposes irrelevant, sign error in 
{\em loc.cit.}
We can apply the antiautomorphism $ * $ to get a similar formula with the
$ T_{(c,m)} $ operating on the right.

\medskip

Let $ \lambda^{i} $ 
be the partition of $ n+1 $ obtained from $ \lambda $ by adjoining a node
to the $ i $'th row (if possible).
Let $t^i$ be the $ \lambda^{i} $-partition obtained from 
$t_{\lambda}$ by filling in the extra node 
with $ n+1 $. 

\medskip

Then one gets from the above formula (\ref{Murphyformula}), 
or rather its $*$-version,
and multiplication with an
appropriate $ T_{w} $ that the 
$$ x_{\lambda \lambda}\, T_{w_{\lambda}} \, y_{(t^{i},s)^{\prime}} $$
will all occur 
in the induced module
for all $  \lambda^{i} $-tableaux $ s $. 

\medskip 
\noindent

We finally claim that 
\begin{equation}\label{basis}
 \{ \, x_{\lambda \lambda}\, T_{w_{\lambda}} \, y_{(t^i,s)^{\prime} } \, \mid s \, 
\,{\rm standard } \, \, \lambda^{i} 
{\rm -tableaux }, \,\, i = 1,\ldots \, \} 
\end{equation}
is a linearly independent set. To do so we may specialize $ q=1 $; then we
have that 
$$ x_{\lambda \lambda}\, T_{w_{\lambda}} \, y_{ (t^i,t^i)^{\prime}} \, = \, 
x_{\lambda \lambda}\, T_{w_{\lambda}} \, y_{\lambda^{\prime},
\lambda^{\prime}}\,\left( \,1 - T_{(j,n)} - \ldots - T_{(m-1,n)}\,\right) $$
and since the $(i,n)$ are coset representatives of $ S_{n-1} $ in
$S_n$ we deduce that at least the 
$$ x_{\lambda \lambda}\, T_{w_{\lambda}} \, y_{ (t^i,t^i)^{\prime}}  $$ are 
linearly independent with $ i $ varying. 
Now letting $s$ vary over standard $\lambda $-tableaux the 
$$  \, x_{\lambda \lambda}\, T_{w_{\lambda}} \, y_{\lambda^{\prime},
s^{\prime}} $$ are linearly independent (a basis of the Specht module),
so the argument used before
can be generalized. 

\medskip

We have now shown that the elements from (\ref{basis}) form 
a basis of the module from the Lemma. Defining
$$ V_i := 
 \{ \, x_{\lambda \lambda}\, T_{w_{\lambda}} \, y_{(t^i,s)^{\prime} } \, \mid s \, 
\,{\rm standard } \, \, \lambda^{k} 
{\rm -tableaux }, \,\, k = 1,\ldots , i \, \} 
$$ 
it follows from Lemma 3.5 of [M1] that the $V_i $'s 
are $ H_{n+1} $ modules, and since clearly $V_i/V_{i-1} \cong 
S(\lambda^i) $, the proof of the Lemma is finished.
\end{pf*}

\medskip 
\noindent

We now come to the diagonalizing procedure; in the rest of this section 
we shall be working over $ {\mathbb Q }(q) $. 
The idea is to reuse Murphy's 
construction of Young's seminormal form, based on the idempotents $ E_t $. 
So let us recall the definition and basic properties of these.

\medskip 
\noindent

Let $ \lambda \vdash n $ and let $ t $ be a $ \lambda $-tableau. 
The (generalized) residue of the $ (i,j) $ node of 
$ t $ is defined to be $ [j-i]_q $ where 
$$ [k]_q := 1 + q + \ldots + q^{k-1} $$
The residue of the node occupied by $ m $ in $ t $ is denoted $ r_t(m) $ and 
the set of possible residues for standard tableaux by 
$ {\cal R}(m) $.
 
\medskip 
\noindent

For any tableau $ t $ we define
$$ E_t \, := \, \prod_{m=1}^n \, \prod_{c \in {\cal R}(n) \setminus \{ r_t(m) \}}
\, \frac{L_m-c}{ r_t(m)-c} $$
where $ L_m $ is the $q$-analogue of the 
Murphy-Jucys operator introduced in [M1]:
$$ L_m = q^{-1} T_{(m-1,m)} + 
q^{-2} T_{(m-2,m)} + \ldots + 
q^{-m} T_{(1,m)} $$
The set $$ \{\, E_t \, \mid \, t \, \, { \rm standard }\, \} $$ is a 
complete set of 
orthogonal idempotents of the Hecke algebra $ H_n  $ 
(defined over $ {\mathbb Q }(q) $) while $ E_t = 0 $ for
$ t $ nonstandard.
 
\medskip 
\noindent
 
Furthermore we 
have the following key formula
\begin{equation}\label{key}
 y_{st} \, E_u \, = \, \delta_{tu^{\prime}} \, y_{st} + \sum_{ \{ \,\sigma , \tau 
\, \mid \,
\sigma \, \tau \, {\rm standard } \,  (\sigma, \tau ) > ( s, t )  \, \} }
\, a_{ \sigma \tau } \,  y_{  \sigma \tau } 
\end{equation}
where $ a_{ \sigma \tau } \in {\mathbb Q }(q) $. 
(This version of the formula is obtained by combining Theorem 4.5 and
the $\eta $-version of (5.1) of [M1]).
It implies that if we 
set $ f_t \, := \, z_{\lambda t }\, E_t $ we get a new basis 
(Young´s seminormal form) of the Specht module 
$$\{ \,f_t \, \mid
\,t \,\, { \rm standard } \, \lambda { \rm -tableau } \} $$

\medskip 
\noindent

The orthogonality of the idempotents implies that the $ f_t $ are orthogonal
with respect to the bilinear form on $ S(\lambda) $.
Moreover one can calculate the length of the $ f_t $; here once again 
one of the main ingredients is the above formula. All this theory is 
developed in the papers [M1],[M2] and [M3]. 

\medskip 
\medskip 

We now focus on the subquotient $ S(\lambda^i) $ of 
$ {\rm Ind}\, S(\lambda) $ which arises 
from adjoining a node to the $ i $'th row of $ \lambda $. Recall the 
$ \lambda^i $-tableau $ t^i $ which is obtained from $ t^{\lambda} $ by 
letting the 
new node be an $ n+1 $ node. Let $ d^i := d(t^i)^{-1} $. 
The lowest vector of the subquotient $ S(\lambda^{i}) $ in 
$ {\rm Ind}\, S(\lambda) $ is now the coset of 
\begin{equation}\label{coset}
 x_{\lambda \lambda}\, T_{w_{\lambda}} \,T_{d(t^i)}^* \, 
y_{(\lambda^i, \lambda^i)^{\prime}} \, = \, 
x_{\lambda \lambda}\, T_{w_{\lambda}} \,T_{d^i} \, 
y_{(\lambda^i, \lambda^i)^{\prime}}
\end{equation}

\medskip 
Let $ t_i $ be the lowest $ \lambda^i $-tableau, in which the 
$ \{1, \ldots , n+1 \} $ are filled in increasingly along the columns.
(Note that $ t_i $ is {\em not} equal to the conjugate of $t^i$).
Then 
\begin{equation}\label{*}
 *\,\,\, x_{\lambda \lambda}\, T_{w_{\lambda}} \,T_{d^i} \, 
y_{(\lambda^i, \lambda^i)^{\prime}}
 \, E_{t_i} 
\end{equation}
is equal to $ x_{\lambda \lambda}\, T_{w_{\lambda}} \,T_{d^i} \, 
y_{(\lambda^i, \lambda^i)^{\prime}} $ modulo the submodule 
$ V_{i-1} $ of 
$ { \rm Ind} \, S(\lambda) $ 
corresponding to the $ \lambda^j $'s with $ j < i $. Furthermore, $ * $ is 
orthogonal to $ V_{i-1} $; here the bilinear form on 
${ \rm Ind} \, S(\lambda) $ is the restriction of the one on the Hecke 
algebra $ H_{n+1} $ -- it is given by 
$$ \langle a,b \rangle := { \rm coefficient \,  of }\,  T_1 \, 
{ \rm in } \, ab^* $$

\medskip 

Let $ U_i $ be the orthogonal complement of $ V_{i-1} $ in $ V_i $. Then 
$ U_i \cong S(\lambda^i) $ canonically and --
like in the restriction functor case -- we get two forms on $ U_i $,
which we wish to compare.

\medskip 

We can view $ * $ as the lowest vector of $ U_i $; it is mapped to  
\begin{equation}\label{**}
 ** \, \, \, x_{\lambda^i \lambda^i}\, T_{w_{\lambda^i}}  \, 
y_{(\lambda^i, \lambda^i)^{\prime}}
\end{equation}
under the above isomorphism $ U_i \cong S(\lambda^i) $. So our task is to 
calculate the lengths of $ * $ and $ ** $ 
and compare. The quotient of the two lengths is the number we are looking 
for. 

\medskip 
\medskip

We need a further piece of notation. 
Let $ \lambda \vdash n $. Denote by $ h_{ij} $ the hook length of the hook
centered at $ (i,k) $ and define the hook product $ h_{\lambda}$ by 
$$  h_{\lambda} = \prod_{(i,j)\in \lambda}\, [h_{ij}]_q $$ 

\medskip
Let $ t $ be a $ \lambda $-tableau; if $ n $ is in the $i$th row of $t $ 
we define the hook-quotient $ \gamma_{tn} $ to be 
$$ \gamma_{tn} = \prod_{j=1}^{\lambda_i}\, \frac{[h_{ij}]_q}{[h_{ij}-1]_q} $$
For $ m<n $ we define $ \gamma_{tm}$ similarly, except that we this time 
remove all the nodes $ m+1, \ldots \, n $ from $ t $. Finally we let 
$$ \gamma_t = \prod_{m=2}^{n}\, \gamma_{tm} $$ 
In the case of $t^{\lambda}$, we then get 
$$ \gamma_{\lambda} = \prod_{i>0}\, \prod_{j=1}^{\lambda_i} \, [j]_q $$

We abbreviate $ \lambda_i^{\prime} $ for $ (\lambda^i)^{\prime} $.
We now have the following Lemma:
\begin{lemma}
Let $ * $ be as above. Then $$ \langle *,* \rangle \, = \, 
\gamma_{\lambda_i^{\prime}}\,
\gamma_{\lambda}\,\prod_{j=1}^{\lambda_i} \, 
\frac{[h_{ij}-1]_q}{[h_{ij}]_q} $$
\end{lemma}
\begin{pf*}{Proof}
First of all we see from Theorem 4.5, including the remarks after the
Theorem, and Lemma 6.1 in [M1] that 
$$ * \, = \, q^{n - \alpha(\lambda)}\,\gamma_{\lambda_i^{\prime}} \,
x_{\lambda \lambda}\, T_{w_{\lambda}} \,T_{d^i} \,  E_{t_i} $$
where $ \alpha $ is the function on partitions (and tableaux) introduced 
in the proof of the preceding Lemma. However, since we are really only 
interested in a certain valuation of $\langle *,* \rangle $, we shall from 
now on omit the $ q $-power of the expression. Hence we get
$$ \begin{array}{ccl}
 \langle *,* \rangle &  \, = \, & \gamma_{\lambda_i^{\prime}}^2 \,
\langle \, x_{\lambda \lambda}\, T_{w_{\lambda}} \,T_{d^i} \,  E_{t_i},\,
x_{\lambda \lambda}\, T_{w_{\lambda}} \,T_{d^i} \,  E_{t_i}\, \rangle \\ & & 
\gamma_{\lambda}\,\gamma_{\lambda_i^{\prime}}^2
\langle \, x_{\lambda \lambda},\, T_{w_{\lambda}} \, T_{d^i}\, E_{t_i}\,
 T_{d^i}^* \,T_{w_{\lambda}}^* \,\,\rangle
\end{array} $$
where we used that $ x_{\lambda \lambda}^2 = \gamma_{\lambda}\, 
x_{\lambda \lambda} $ which follows from the definitions (or from Theorem 
4.5 and Lemma 6.1 in [M1]). 

\medskip
In order to calculate $ T_{w_{\lambda}} \, T_{d^i}\, E_{t_i}\,
T_{d^i}^{\,*} \,T_{w_{\lambda}}^{\,* }$ we shall use Theorem 6.4 of [M1], 
which gives 
formulas for the product $ \zeta_{us}\, T_v $ where $ v = (i-1,i) $ and 
$ \zeta_{us} := E_u \, x_{us} \, E_s $ (actually [M1]
uses the basis $ \xi_{us} $ in the definition of $ \zeta_{us} $; however 
this basis is related to the standard 
basis $ x_{us} $ by an upper triangular matrix with ones on the diagonal and 
thus gives rise to the same $ \zeta_{us} $). To get the action 
on the $ E_t $'s
we then use the formula 
$$ \zeta_{tt} = \gamma_{t} \, E_t $$
which is also proved in [M1].

\medskip
Let us quote the mentioned Theorem 6.4.
\begin{thm} Let $s, \, \, u $ be standard tableaux of the same shape, $ v = 
(i-1,i ) $, $ t= sv $. Let $ (a,b) $ and $ ( a^{\prime}, b^{\prime} ) $ be 
the nodes occupied by $ i-1 $ and $ i $ respectively in $ s $ and let 
$ h = b - b^{\prime} - a + a^{\prime} $; then  \medskip
$$ \zeta_{us} \, T_v  = \left\{ \begin{array}{ll} 
      \frac{1}{[h]_q}\, \zeta_{us} & if \, \, \, | h | = 1,  \\ \\
      \frac{1}{[h]_q}\,\zeta_{us} \,+\, \zeta_{ut} &  if \, \, \, h>1, \\ \\
      \frac{1}{[h]_q}\,\zeta_{us}  \, + \, \frac{q\, [h+1]_q  \, [h-1]_q }
      {[h]^2_q}\, \zeta_{ut} & if \, \, \, h<-1, 
                              \end{array}
                 \right. $$
\end{thm}

By applying $ * $ we get formulas for the left multiplication with $ T_v $
as well. 

\medskip 

The plan is now to write $ d_i $ and $ w_{\lambda} $ as a product of basic 
transpositions and use the above formulas.
Let $s^i$ be the $ \lambda^i $-tableau obtained from $ t^{\lambda} $ 
by letting the new 
node be an $ n+1 $-node. 
Recalling that $ t_i $ is 
the smallest $ \lambda^i $--tableau, the result of applying this to 
$ T_{w_{\lambda}} \, T_{d^i}\, E_{t_i}\,
T_{d^i}^{\,*} \,T_{w_{\lambda}}^{\,* }$ will be a 
sum of $ \zeta_{us} $ with $ u $ and $ s $ less than $ s^i $.

\medskip 

However, since 
\begin{equation}\label{ortho} x_{\lambda \lambda} \,
E_{\tau} = E_{\tau}\,x_{\lambda \lambda} = 
0 \,\,\,\,{\rm for }  \,\,\tau < \lambda  
\end{equation}
(by Theorem 4.5, Lemma 6.1 and (5.5) of [M1]) and since 
$ \zeta_{us} := E_u \, x_{us} \, E_s $ we actually only need a small part of 
the sum to determine $ \langle *,* \rangle  $.

\medskip

Notice first of all that 
$$ \begin{array}{lll}
d^i & = & ( k,k+1, \ldots, n+1 ) \\
    & = & ( n+1,n)\, (n,n-1) \, \cdots \, ( k+1,k ) $$
   \end{array} $$
where $k$ is the number in the position $ ( i, \lambda_i +1) $
of $ t^{\lambda^i} $. This is a reduced presentation 
of $ d^i$ and so $ T_{d^i} = T_{(n+1,n)}\, \cdots \,T_{(k,k+1)} $. 
The product 
$ T_{d^i}\, E_{t_i}\, T_{d^i}^* $ 
will therefore only involve $ \zeta_{us} $'s
in which $ n+1 $ occurs in $ u $ and $ s $ in positions higher than in $ s^i $;
removing $ n+1 $ from these will lead to partitions smaller than or 
equal to $ \lambda $. But by (\ref{ortho}) we can neglect those
$ u $ and $ s $ where $ n+1 $ is in a strictly higher position. 

\medskip
We can repeat this argument on the remaining numbers $ 1,2, \ldots, n $ and 
find that only for $ u=s=s^i $ there will be a contribution 
to $ \langle * , * \rangle $ from
$ T_{w_{\lambda}} \, T_{d^i}\, E_{t_i}\,
T_{d^i}^{\,*} \,T_{w_{\lambda}}^{\,* }$ 

\medskip
Let $ s $ and $ t $ be standard tableaux with $ t = s \, v, \, \, \, s<t $ for 
$ v= ( i, i-1) $. Then we are in the situation $ h < \, -1 $ of the Theorem,
so we get (once again ignoring the $ q $-power):
$$ \begin{array}{lll}
\frac{1}{\gamma_s}\,T_v\,E_s\,T_v & = & \frac{1}{\gamma_s^2}\,T_v\,
\zeta_{ss}\,T_v \\
&= & \frac{1}{\gamma_s^2}\,\left(\frac{[h+1]_q\,[h-1]_q}
{[h]_q^2} \right)^2
\, \zeta_{tt} \, + {\rm lower \, \,  terms} \\
& =& \frac{1}{\gamma_t^2}\,\zeta_{tt}\, + {\rm lower \, \,  terms} \\
& =& \frac{1}{\gamma_t}\, E_t + {\rm lower \, \,  terms} 
\end{array}$$

\medskip

We are now in position to calculate the length $ \langle *, * \rangle $. 
Using the above formula we get 
$$ \begin{array}{lll}
\langle *,* \rangle &  \, = \, & 
\gamma_{\lambda}\,\gamma_{\lambda_i^{ \prime}}^2
\langle \, x_{\lambda \lambda},\, T_{w_{\lambda}} \, T_{d^i}\, E_{t_i}\,
 T_{d^i}^* \,T_{w_{\lambda}}^* \,\,\rangle \\
& = & \frac{\gamma_{t_i}}{\gamma_{s^i}}\,\gamma_{\lambda}\,
\gamma_{\lambda_i^{ \prime}}^2 \, \langle \, x_{\lambda \lambda}, E_{s^i}
\rangle \\
& = & \frac{\gamma_{t_i}}{\gamma_{s^i}}\,{\gamma_{\lambda}}^2\,
\gamma_{\lambda_i^{\prime}}^2 \, \langle  E_{s^i}, E_{s^i} \rangle
\end{array} $$
For the last equality we factored $ E_{\lambda} $ out of $ E_{s^i} $,
applied the formula $ x_{\lambda \lambda} \, E_{\lambda} = 
\gamma_{\lambda}\, x_{\lambda \lambda} $ and plugged in the extra factor 
once again.

\medskip
But the length of the idempotents is known, see the last page of [M1], 
it is 
$$ \langle E_t,E_t \rangle = \frac{1}{h_{\lambda}} = \frac{1}{\gamma_t \,
\gamma_{t^{\prime}}} $$
This holds at least for $ t = t^{\lambda} $ by {\em loc.cit.} However one 
checks 
that the argument of Theorem 6.6 actually is valid for any $\lambda$-tableau
$ t $ and thus the length $ \langle E_t,E_t \rangle $ is independent of $t$.

\medskip
 
All in all we have
$$ \begin{array}{lll}
\langle *,* \rangle &  =  &
\frac{\gamma_{t_i}}{\gamma_{s^i}}\,{\gamma_{\lambda}}^2 \,
{\gamma_{\lambda_i^{\prime}}}^2 \frac{1}{\gamma_{\lambda_i^{\prime}}\,\gamma_{t_i }} \\
& = & \gamma_{\lambda_i^{\prime}}\,
\gamma_{\lambda}\,\prod_{j=1}^{\lambda_i} \, 
\frac{[h_{ij}-1]_q}{[h_{ij}]_q}
\end{array} $$
and the proof of the Lemma is finished.
\end{pf*}

Let us return to $ {\rm Ind} \, S(\lambda ) $ and its filtration 
$$ 0 \, \subseteq \, V_1 \, \subseteq V_2 \, \subseteq \, \ldots \,  V_m \, = 
\, I 
\,\, {\rm{with}} \,\,  V_i \, /\, 
V_{i-1} 
\cong S(\lambda^i ) $$ 
The projection map goes as follows:
$$ \begin{array}{l}
\pi_{i}: \, \, V_i \, \rightarrow \, \, S(\lambda^i) \\
x_{\lambda \lambda}\, T_{w_{\lambda}} \, y_{ (t^j t^j)^{\prime} } \, \, 
\mapsto \left\{
\begin{array}{ll}
x_{\lambda^i \lambda^i}\, T_{w_{\lambda^i}} \, y_{(\lambda^i, 
\lambda^i)^{\prime} } &   if \,\,j=i \\
0                      &   if \,\, j>i 
\end{array} \right.
\end{array} $$ 
We already saw that 
$ x_{\lambda \lambda}^2 = \gamma_{\lambda}\,x_{\lambda \lambda} $ and
$ y_{\lambda \lambda}^2 = \gamma_{\lambda}\,y_{\lambda \lambda} $. Since 
furthermore $ x_{\lambda \lambda}\,y_{t_{\lambda} t_{\lambda}}\, = \, 1 $ --
ignoring $ q $-powers -- we find
\begin{equation}\label{1}
 \langle \, x_{\lambda^i \lambda^i}\, T_{w_{\lambda^i}} \, 
y_{(\lambda^i, \lambda^i)^{\prime} }, \,
x_{\lambda^i \lambda^i}\, T_{w_{\lambda^i}} \, 
y_{(\lambda^i , \lambda^i)^{\prime} }
\, \rangle = \gamma_{\lambda^i} \,
\gamma_{(\lambda^i)^{\prime}} 
\end{equation}
A similar argument shows, using (\ref{coset}), that 
\begin{equation}\label{2}
 \langle
x_{\lambda \lambda}\, T_{w_{\lambda}} \, y_{ t^i t^i }\,, 
x_{\lambda \lambda}\, T_{w_{\lambda}} \, y_{ t^i t^i } \rangle = 
\gamma_{\lambda} \,\gamma_{(\lambda^i)^{\prime}} 
\end{equation}
So we conclude that the isomorphism $ V_i \, /\, V_{i-1} \cong S(\lambda^i )$
stretches all squared lengths by the factor $ [\lambda_i +1]_q $.

\medskip
Recall $ U_i $ the orthogonal complement of $ V_{i-1} $ in $ V_i $, we have 
$ U_i \cong S(\lambda^i) $. We consider then two $S_{n+1}$-invariant forms on 
$U_i$, namely the usual one $ \langle *,* \rangle $ coming from the embedding
$ U_i \subseteq {\rm Ind}\,S(\lambda) $ and $ \langle *,* \rangle_1 $ 
which is the pullback of the form on $ S(\lambda^i) $ along $ \pi_i $.
The two forms differ by a constant $ m $: $$ \langle *,* \rangle =
m \, \langle *,* \rangle_1 $$
\medskip
The next theorem gives the promised statements about $ m $: 

\begin{thm} Let $l>1$ and let $\nu_l$ be the $ [l]_q $-adic valuation 
$ { \mathbb Z} [q,q^{-1}] $ defined in (\ref{valuation}). 
Consider the $l$-diagram of $ \lambda^i $ and 
assume that $ \tau = \lambda^i \setminus \lambda $ is a $k$-node.
Then $$\nu_q(m)\,= N^l_i(\lambda^i,\lambda) $$
where $ N^l_i(\lambda^i,\lambda ) $ is the number of removable $k$-nodes 
situated to the left of $ \tau $ minus the number of indent $ k $-nodes
situated to the left of $ \tau $. 
\end{thm}

\begin{pf*}{Proof}
This follows 
by combining the previous Theorem 
with (\ref{*}), (\ref{**}), (\ref{1}), (\ref{2}), 
along with an argument 
similar to the one of Theorem 1.  
\end{pf*}

\medskip 

\section{\bf Duality}
We now return to the symmetric group case. We wish to investigate the 
above logic on the duality of the symmetric group. 

\medskip 
Recall that for any right module $ M $ of a finite group $ G $ the dual module 
is defined by 
$$ M^*:={\rm Hom}_k\, (M,k);\, \, (fg)(m)=f(mg), \, \, f \in M^*, \, 
m \in M, \, g \in G $$
There is then the following classical result:
$$ S(\lambda)^* \cong S(\lambda^{\prime}) \otimes k_{alt} $$

The way James [J] proves this result is as follows: $ S(\lambda) $ is by 
construction a submodule of the permutation module $ M(\lambda ) $, which 
comes with an invariant form $ \langle *,* \rangle_{\lambda} $, 
which is nondegenerate 
independently of the field $ k $. One can construct a surjection 
$$ \pi: \, \, M(\lambda) \twoheadrightarrow S(\lambda^{\prime}) 
\otimes k_{alt} $$ 
and can prove that $ {\rm ker}\,\pi = S(\lambda)^{\bot}$. 
Since $ M(\lambda)/S(\lambda)^{\bot} \cong S(\lambda)^* $, 
the isomorphism 
$ S(\lambda)^* \cong S(\lambda^{\prime}) \otimes k_{alt} $ now follows.

\medskip

Working over $ \mathbb Q $ the restriction of $ \pi $ to $ S(\lambda ) $ is 
an isomorphism, i.e. $ S(\lambda)^* \cong  S(\lambda ) $. So in that case 
there are two invariant 
forms on $ S(\lambda^{\prime}) \otimes k_{alt}$: the first one is 
$ \langle *, \rangle_{\lambda^{\prime}} \otimes 1 $ which comes from the 
embedding $ S(\lambda^{\prime}) \subseteq M(\lambda^{\prime}) $, the second 
one is obtained by carrying over the form on $ S(\lambda) $ to 
$ S(\lambda^{\prime}) \otimes k_{alt} $.
As usual the two forms differ by a scalar $ m \in \mathbb Q $; 
we wish to calculate the $ \nu_{q} $-adic valuation of it. 

\medskip

Let us therefore take a closer look at the way $ \pi $ is defined. For a 
$ \lambda $-tableau $ t $ we denote by $\kappa_t $ the alternating column 
sum $ \kappa_{t}= \sum_{C(t)}\,{\rm sgn}\sigma \, \sigma $ and by $ \rho_{t} $
the row sum $ \rho_{t}= \sum_{R(t)}\,\sigma $ of $ t $. So our generators 
$ e_t $ of the Specht module $ S(\lambda) $ have the form
$ e_t = \{t\}\, \kappa_t $ -- here $ \{t\} $ was the tabloid class of $ t $.

\medskip

With this notation, $ \pi $ is the map given by 
$$ \begin{array}{l}
\pi: \, S(\lambda) \longrightarrow S( \lambda^{\prime} ) \otimes  k_{alt} \\
 e_t = \{t\} \mapsto \{t^{\prime}\}\,  \kappa_t^{\prime} \, \rho_t^{\prime}
\, \otimes \,1
\end{array} $$

Using $\rho_{t}^2 = \vert R(t) \vert \, \rho_{t} $, we find that 
$$  \begin{array}{lll}
\langle \, \pi(e_t), \pi(e_t)\, \rangle & = & \langle \,
 \{t^{\prime}\}\,  \kappa_t^{\prime} \, \rho_t^{\prime}, 
 \{t^{\prime}\}\,  \kappa_t^{\prime} \, \rho_t^{\prime}\, \rangle \\
& = & \langle t^{\prime}, t^{\prime}\kappa_t^{\prime} \, {\rho_t^{\prime}}^2\,
\kappa_t^{\prime}\, \rangle \\
& = & \vert R(t)\vert \, \langle \, t^{\prime}, t^{\prime}\kappa_t^{\prime} \, \rho_t^{\prime}\,
\kappa_t^{\prime}\,\, \rangle \\
& = & \vert R(t)\vert \,h_{\lambda^{\prime}}\, \langle\, t^{\prime}, 
t^{\prime} 
\kappa_t^{\prime}\,\, \rangle \\
& = & \vert R(t)\vert \,h_{\lambda^{\prime}}\\
\end{array} $$
One again $h_{\lambda} $ denotes the product of all hook lengths in 
$ \lambda $; for the second last equality we used Lemma 23.2 in James's book
[J]. 

\medskip 

On the other hand we have that:
$$ \langle \, e_t, e_t \, \rangle = \langle \, 
\{t\} \kappa_t, \{t\}\kappa_t \, \rangle = \vert C_t \vert $$
Since $$ \vert C_t \vert = \vert R_{t^{\prime}}  \vert $$ we deduce that 
our constant $ m $ is equal to $ h_{\lambda} $. 
But by 2.7.40 the number of $p$-divisible hooks 
equals the $p$-weight of $ \lambda $, i.e.
$$ \nu_q(m) = \vert \lambda \vert_p $$
But this is exactly the $q$-power that appears in [LLT]'s duality, see 
Theorem 7.2.of {\em loc. cit.}

\section{\bf The ${SL}_2 $-situation.}

In this section we shall comment on the problem of determining the tilting 
modules outside the lowest $ p^2 $-region. This is probably the most 
difficult problem in the modular representation theory, there is so far 
not even a conjecture around. In type A this problem is known to be 
equivalent to determining the decomposition matrix of the symmetric group.

\medskip

In fact, if $ T(\lambda) \, \in \, \mbox{mod-}Gl_m  $ is the tilting module 
labeled by the partition $ \lambda $ of $n $,
$ m \geq n$, we have, see [E]:
\begin{equation} [T(\lambda), \Delta(\mu) ] = [S(\mu), D(\lambda) ] 
\, \, \, \, \, \, \, \lambda \in \mbox{Par}_p, \mu \in \mbox{Par} 
\end{equation}
where $\mbox{Par}$ denotes the partitions of $n$ while 
$\mbox{Par}_p$ denotes the $p$-regular partitions of $n$.
Based on this formula, or rather a $q$-analogue of it 
relating the relevant $q$-Schur and Hecke-algebras,
one can show that the LLT-algorithm is 
nothing but an extension of Soergel's tilting algorithm 
to singular weights.

\medskip

Now the calculations done in the previous chapters lead naturally to the 
following guess on how to obtain the general tilting modules: {\it perform the 
LLT-algorithm as usual but replace everywhere the quantum valuations by
ordinary $p$-adic valuations. }

\medskip

In this section we shall see that this 
idea seems to work
at least in the $ {SL}_2 $-situation.
On the other 
hand we point out right from the beginning that already in the 
$ {SL}_3 $-case the direct generalization of the 
resulting algorithm does {\em not} work. See however the work of 
A. Cox, [AC].

\medskip
\medskip

Let us recall the work of K. Erdmann [E] on the tilting modules for 
$ {SL}_2 $. We start out with some notation: any natural number
$ a $ has a $p$-adic decomposition 
$$ a = a_0 + a_1 p +  a_2 p^2 + \ldots a_k p^k $$
For any other natural number $ b $, we shall say that $ a $ contains $b$ 
if in the $p$-adic decomposition of $b $ 
$$ b= b_0 + b_1 p +  b_2 p^2 + \ldots b_l p^l $$
$l<k $ and for all $ i $: $b_i = a_i $ or $ b_i = 0 $. 

\medskip

The representation theory of $ {SL}_2 $ is parameterized by 
$ \mathbb N $, i.e. for every natural number $n$ there is a tilting module 
$T(n)$ with highest weight $n$ and so on. K. Erdmann has now obtained 
the following result, [E]:

\begin{thm}Let $ p \not=  2 $. Then the multiplicity of the Weyl module 
$ \Delta(s) $ in $ T(m) $ is given by the formula:
$$ [ T(m),\Delta(s)] = \left\{ 
\begin{array}{ll}
1 & \mbox{ if } m+1 \mbox{ contains }  \frac{1}{2}(m-s) \\
0 & \mbox{ otherwise } 
\end{array}
\right. $$
\end{thm}

\noindent
{\bf Example}: Let us take $ p=3 $ and $ m= (p^3 + p^2 +2 ) -1 $. The 
Theorem then gives rise to the following alcove pattern, where the dots indicate 
Weyl composition factors in $ T(m) $ -- all multiplicities are one.

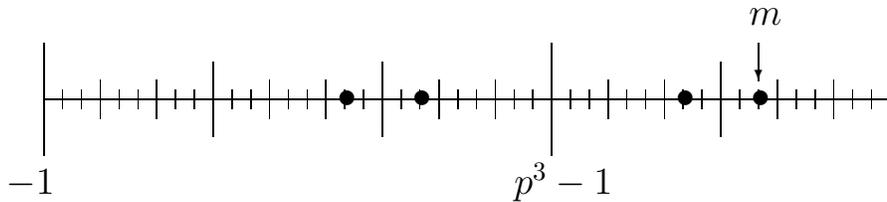
\setlength{\unitlength}{0.25cm}
\begin{center}
\begin{picture}(5,5)
\put(-20,0){\line(1,0){45}}
\linethickness{0.075mm}
  
\multiput(-20,-0.5)(1,0){46}{\line(0,1){1}}
  \multiput(-20,-1)(3,0){16}{\line(0,1){2}} 
  \multiput(-20,-2)(9,0){6}{\line(0,1){4}}
  \multiput(-20,-3)(27,0){2}{\line(0,1){6}} 
  \put(-22,-5){$-1$}
  \put(5,-5){$p^3-1$}
  \put(17.6,-0.4){$\bullet$}
  \put(13.6,-0.4){$\bullet$}
  \put(-0.4,-0.4){$\bullet$}
  \put(-4.4,-0.4){$\bullet$}
  \put(17.5,4){$m$}
  \put(18,3){\vector(0,-1){2}}
\end{picture}
\end{center}

\medskip
\medskip
\medskip
\medskip
\medskip
\medskip
\medskip
\medskip
\medskip
The rule for obtaining the picture is the following: reflect first $ m $ in 
the last $ p $-wall before $ m $, then in the last $ p^2 $-wall before 
$ m $ and so on. Finally ``symmetrize'' the picture. 
The number of 
$\Delta$-factors in $ T(m ) $ hence equals $ 2^l $, where $ l = \#
\{ i | a_i \not= 0 \} $ for $ m+1 = a_0 + a_1\,p + \ldots +  a_k\,p^k $.

\medskip
We shall present the effect of the LLT-algorithm 
in a series of examples.
The representation theory of $ Gl_2 $ is parameterized by Young diagrams 
with at most two lines. The diagram 
$ (\lambda_1, \lambda_2 )$
passes under the restriction to the $  Sl_2 $-weight $ m=\lambda_1 - \lambda_2 $.
Let us therefore consider the LLT algorithm on Young-diagrams with at most two lines.
It adds nodes in different ways and we simply neglect 
all diagrams with more than two lines appearing. It only involves 
the $ f_i $-operators, which when adding a node to the second line may
produce a $ q $-power (namely $q^{N^r_i}$, where  ${N^r_i}$ is $ 1 $ if the first 
line has an addable $ i $-node, $ -1 $ if it has a removable $ i $-node).

\medskip
\noindent
{\bf Example 1}: Take $ l= 5 $, $ (\lambda_1,\lambda_2) = ( 15,9 ) $. The 
corresponding $ Sl_2 $ weight is $ m= 6 $. 

\medskip

\setlength{\unitlength}{0.5cm}

\begin{picture}(3,3)

\multiput(5,0)(1,0){1}{\framebox(0.9,0.9){$0$}}
\multiput(6,0)(1,0){1}{\framebox(0.9,0.9){$1$}}
\multiput(7,0)(1,0){1}{\framebox(0.9,0.9){$2$}}
\multiput(8,0)(1,0){1}{\framebox(0.9,0.9){$3$}}
\multiput(9,0)(1,0){1}{\framebox(0.9,0.9){$4$}}
\multiput(10,0)(1,0){1}{\framebox(0.9,0.9){$0$}}
\multiput(11,0)(1,0){1}{\framebox(0.9,0.9){$1$}}
\multiput(12,0)(1,0){1}{\framebox(0.9,0.9){$2$}}
\multiput(13,0)(1,0){1}{\framebox(0.9,0.9){$3$}}
\multiput(14,0)(1,0){1}{\framebox(0.9,0.9){$4$}}
\multiput(15,0)(1,0){1}{\framebox(0.9,0.9){$0$}}
\multiput(16,0)(1,0){1}{\framebox(0.9,0.9){$1$}}
\multiput(17,0)(1,0){1}{\framebox(0.9,0.9){$2$}}
\multiput(18,0)(1,0){1}{\framebox(0.9,0.9){$3$}}
\multiput(19,0)(1,0){1}{\framebox(0.9,0.9){$4$}}

\multiput(5,-1)(1,0){1}{\framebox(0.9,0.9){$4$}}
\multiput(6,-1)(1,0){1}{\framebox(0.9,0.9){$0$}}
\multiput(7,-1)(1,0){1}{\framebox(0.9,0.9){$1$}}
\multiput(8,-1)(1,0){1}{\framebox(0.9,0.9){$2$}}
\multiput(9,-1)(1,0){1}{\framebox(0.9,0.9){$3$}}
\multiput(10,-1)(1,0){1}{\framebox(0.9,0.9){$4$}}
\multiput(11,-1)(1,0){1}{\framebox(0.9,0.9){$0$}}
\multiput(12,-1)(1,0){1}{\framebox(0.9,0.9){$1$}}
\multiput(13,-1)(1,0){1}{\framebox(0.9,0.9){$2$}}
\end{picture}

\medskip
\medskip
\medskip
\medskip
\medskip
\medskip
\medskip
The last two nodes of the $ 5 $-diagram are different, corresponding to 
$ m $ being a regular $ Sl_2 $-weight. We can add a $ 0 $-node and a 
$ 3 $-node, i.e. only $ f_0 $ and $f_3 $ operate non trivially on the 
diagram. The action is as follows:

$$ f_0\, ( \lambda_1, \lambda_2 ) = ( \lambda_1 +1, \lambda_2 ),
\,\,\, { \rm i.e}
\,\,\,\,\,  f_0\, ( m ) = (m+1) $$
$$ f_3\, ( \lambda_1, \lambda_2 ) = ( \lambda_1 , \lambda_2 +1 ),
\,\,\, { \rm i.e}
\,\,\,\,\,  f_3\, ( m ) = (m-1) $$

\medskip
\medskip

\noindent
{\bf Example 2}: Let $ (\lambda_1,\lambda_2 ) $ satisfy 
$ \lambda_1 - \lambda_2 \equiv -1 \, {\rm mod } \,\, l $, i.e. the 
corresponding $ Sl_2 $-weight is a Steinberg weight. Then the 
$ l $-diagram has the form 

\setlength{\unitlength}{0.5cm}

\begin{picture}(3,3)

\multiput(5,0)(1,0){1}{\framebox(0.9,0.9){$0$}}
\multiput(6,0)(1,0){1}{\framebox(0.9,0.9){$1$}}
\multiput(7,0)(1,0){1}{\framebox(0.9,0.9){$2$}}
\multiput(8,0)(1,0){1}{\framebox(0.9,0.9){$3$}}
\multiput(9,0)(1,0){1}{\framebox(0.9,0.9)}
\multiput(10,0)(1,0){1}{\framebox(0.9,0.9)}
\multiput(11,0)(1,0){1}{\framebox(0.9,0.9)}
\multiput(12,0)(1,0){1}{\framebox(0.9,0.9)}
\multiput(13,0)(1,0){1}{\framebox(0.9,0.9)}
\multiput(14,0)(1,0){1}{\framebox(0.9,0.9)}
\multiput(15,0)(1,0){1}{\framebox(0.9,0.9)}
\multiput(16,0)(1,0){1}{\framebox(0.9,0.9)}
\multiput(17,0)(1,0){1}{\framebox(0.9,0.9)}
\multiput(18,0)(1,0){1}{\framebox(0.9,0.9)}
\multiput(19,0)(1,0){1}{\framebox(0.9,0.9){$r$}}

\multiput(5,-1)(1,0){1}{\framebox(0.9,0.9){$-1$}}
\multiput(6,-1)(1,0){1}{\framebox(0.9,0.9){$0$}}
\multiput(7,-1)(1,0){1}{\framebox(0.9,0.9){$1$}}
\multiput(8,-1)(1,0){1}{\framebox(0.9,0.9){$2$}}
\multiput(9,-1)(1,0){1}{\framebox(0.9,0.9)}
\multiput(10,-1)(1,0){1}{\framebox(0.9,0.9)}
\multiput(11,-1)(1,0){1}{\framebox(0.9,0.9)}
\multiput(12,-1)(1,0){1}{\framebox(0.9,0.9)}
\multiput(13,-1)(1,0){1}{\framebox(0.9,0.9){$r$}}
\end{picture}

\medskip
\medskip
\medskip
\medskip
\medskip
\medskip
\medskip

\noindent
where $ r \equiv  \lambda_1 -1 \equiv \lambda_2 -2 \,\,\,{\rm mod }\,\,
p $. We can thus only 
operate non trivially with $ f_{r+1} $ and obtain the following formula:
$$ f_{r+1}\, ( \lambda_1 ,\lambda_2 ) = ( \lambda_1 +1, \lambda_2 ) +
q\,( \lambda_1 , \lambda_2 +1)    \,\,\,\,\,\, {\rm i.e. }\, \, $$
$$ f_{r+1}\, ( m ) = (m+1) + q\,(m-1)     \,\,\,\,\,\,  \,\,\,\,\,\, 
\,\,\,\,\,\,  \,\,\,\,\,\, \,\\,\,\,\,\,\,$$

\medskip
\medskip
\medskip

\noindent
{\bf Example 3}: Let $ (\lambda_1,\lambda_2 ) $ satisfy 
$ \lambda_1 - \lambda_2 \equiv -2 \,\, {\rm mod } \,\, p $. Then 
$ m+1 $ is a Steinberg weight. We can operate non trivially only with 
$ f_{\lambda_1} $ and find that 

$$ f_{\lambda_1}\, ( \lambda_1, \lambda_2 ) = ( \lambda_1 +1, \lambda_2 )
\,\,\,\,\,\,
{ \rm i.e } \,\,\,\, f_{\lambda_1}\, ( m ) = ( m+1 ) $$

\medskip
\medskip
\medskip
\noindent
{\bf Example 4}: Let $ (\lambda_1,\lambda_2 ) $ satisfy 
$ \lambda_1 - \lambda_2 \equiv 0 \,\, {\rm mod } \,\, p $. Then 
$ m-1 $ is a Steinberg weight and we now get 

$$ f_{\lambda_2}\, ( \lambda_1, \lambda_2 ) = q^{-1}\,
( \lambda_1 , \lambda_2 +1 )\,\,
{ \rm i.e } \,\,\,\, f_{\lambda_2}\, ( m ) =q^{-1}\, ( m-1 ) $$

\medskip
\medskip
\medskip
\noindent
{\bf Example 5}: Let  $ (\lambda_1,\lambda_2 ) $ be as in example 3. We can
then compose the operators of example 3 and 2 to obtain a kind of 
translation through the wall. The result is 

$$ f_{\lambda_1+1}\,f_{\lambda_1}\,( \lambda_1, \lambda_2 ) = 
( \lambda_1 +2 , \lambda_2)+ q \,( \lambda_1 +1 , \lambda_2 +1)  \,\,\,
{ \rm i.e } $$
$$f_{\lambda_1+1}\,f_{\lambda_1}\,(m) = ( m+2 ) +q \,(m)
\\,\,\,\,\,\,\,\,\,\,\,\,\,\,\,\,\,\,\,\,\,\,\,\,\,\,\,\,\,\,\,\,\,\,\,\,\,\,\,\,\,\,\,\,\,\,\,\,\,\,\,\,\,\,\,\,\,$$

\medskip
\medskip
\medskip
\noindent
{\bf Example 6}: Let  $ (\lambda_1,\lambda_2 ) $ be as in example 4. We can
then compose with example 2 and this time translate backwards 
through the wall. The result is 

$$ f_{\lambda_1}\,f_{\lambda_2}\,( \lambda_1, \lambda_2 ) = 
q^{-1}\,( \lambda_1 +1 , \lambda_2 +1 )+ ( \lambda_1 , \lambda_2 +2)  \,\,\,
{ \rm i.e } $$
$$ (m) \mapsto q^{-1}\,(m) +(m-2) \,\,\,\,\,\,\,\,\,\,\,\,\,\,\,\,\,\,\,\,\,\,\,\,\,\,\,\,\,\,\,\,\,\,\,\,\,\,\,\,\,\,\,\,\,\,\,\,\,\,\,\,\,\,\,\,\,\,\,\,\,\,\,\,\,\,\,\,\,\,\,$$

\medskip
\medskip
\medskip

(The combination of these examples proves that the LLT algorithm is equivalent
to Soergel's tilting algorithm in the $ Sl_2 $-case.)

\medskip
Let us now consider the modification of the algorithm mentioned in the
beginning of the section. 
Let $ (\lambda_1, \lambda_2 ) $ be a two line partition with Young 
diagram as follows:

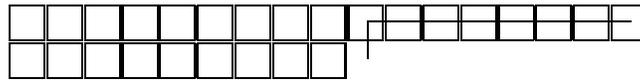
\setlength{\unitlength}{0.5cm}

\begin{picture}(3,3)

\multiput(5,0)(1,0){1}{\framebox(0.9,0.9)}
\multiput(6,0)(1,0){1}{\framebox(0.9,0.9)}
\multiput(7,0)(1,0){1}{\framebox(0.9,0.9)}
\multiput(8,0)(1,0){1}{\framebox(0.9,0.9)}
\multiput(9,0)(1,0){1}{\framebox(0.9,0.9)}
\multiput(10,0)(1,0){1}{\framebox(0.9,0.9)}
\multiput(11,0)(1,0){1}{\framebox(0.9,0.9)}
\multiput(12,0)(1,0){1}{\framebox(0.9,0.9)}
\multiput(13,0)(1,0){1}{\framebox(0.9,0.9)}
\multiput(14,0)(1,0){1}{\framebox(0.9,0.9)}
\multiput(15,0)(1,0){1}{\framebox(0.9,0.9)}
\multiput(16,0)(1,0){1}{\framebox(0.9,0.9)}
\multiput(17,0)(1,0){1}{\framebox(0.9,0.9)}
\multiput(18,0)(1,0){1}{\framebox(0.9,0.9)}
\multiput(19,0)(1,0){1}{\framebox(0.9,0.9)}
\multiput(20,0)(1,0){1}{\framebox(0.9,0.9)}
\multiput(21,0)(1,0){1}{\framebox(0.9,0.9)}

\put(14.5, 0.5){\line(0,-1){1}}
\put(14.5, 0.5){\line(1,0){7}}

\multiput(5,-1)(1,0){1}{\framebox(0.9,0.9)}
\multiput(6,-1)(1,0){1}{\framebox(0.9,0.9)}
\multiput(7,-1)(1,0){1}{\framebox(0.9,0.9)}
\multiput(8,-1)(1,0){1}{\framebox(0.9,0.9)}
\multiput(9,-1)(1,0){1}{\framebox(0.9,0.9)}
\multiput(10,-1)(1,0){1}{\framebox(0.9,0.9)}
\multiput(11,-1)(1,0){1}{\framebox(0.9,0.9)}
\multiput(12,-1)(1,0){1}{\framebox(0.9,0.9)}
\multiput(13,-1)(1,0){1}{\framebox(0.9,0.9)}

\put(14.5, 0.5){\line(0,-1){1}}

\end{picture}

\medskip
\medskip
\medskip
\medskip
\medskip

Let $ h $ be the length of the hook indicated in the diagram. When adding 
a node to the second line, we thus multiply by the $ q$-power with 
exponent $ \nu( \frac{h}{h-1}) $ where $ \nu $ is ordinary 
$ p $-adic valuation.

One now goes through the various examples. Example 5, which was translation 
through the wall, takes the following form:

\setlength{\unitlength}{0.25cm}
\begin{center}
\begin{picture}(5,5)
\put(-25,0){\line(1,0){20}}
\linethickness{0.075mm}
\put(-15,2){\line(0,-1){4}}
\put(-17,-3){$p^a -1$}
\put(-19,1){$1$}

\put(1,-0.3){$\longrightarrow$}

\put(10,0){\line(1,0){20}}
\linethickness{0.075mm}
\put(20,2){\line(0,-1){4}}
\put(18,-3){$p^a -1$}
\put(16,1){$q^a$}
\put(24,1){$1$}

\end{picture}
\end{center}

\medskip
\medskip
\medskip
\medskip
\medskip
\medskip
\medskip
\medskip
\medskip
\noindent
while example 6 takes the form 

\setlength{\unitlength}{0.25cm}
\begin{center}
\begin{picture}(5,5)
\put(-25,0){\line(1,0){20}}
\linethickness{0.075mm}
\put(-15,2){\line(0,-1){4}}
\put(-17,-3){$p^a -1$}
\put(-13,1){$1$}

\put(1,-0.3){$\longrightarrow$}

\put(10,0){\line(1,0){20}}
\linethickness{0.075mm}
\put(20,2){\line(0,-1){4}}
\put(18,-3){$p^a -1$}
\put(16,1){$1$}
\put(24,1){$q^{-a}$}

\end{picture}
\end{center}

\medskip
\medskip
\medskip
\medskip
\medskip
\medskip
\medskip
\medskip
\medskip

We can now start the modified algorithm. There is however one more detail to point 
out. When running the modified algorithm, there will soon be negative 
$ q$-powers occurring. This is not the case in Soergel's algorithm, but 
certainly the case in the LLT-algorithm, generally. Hence, we
shall proceed as LLT do, not only subtracting inductively known tilting
characters when a $ 1 $ is appearing, but also when a negative $ q $-power
is appearing. To be precise: if some negative $ q $-power appears, 
it is possible to subtract an expression on the form $ \gamma \, T $ 
where $ \gamma \in {\mathbb Z }\,[q,q^{-1} ]$ with $ \gamma(q^{-1}) = 
\gamma(q) $ and where $ T $ is a known tilting character to finally 
arrive at an expression involving only positive $ q $-powers.

\medskip
\medskip

Let us work out the case $ p = 3 $ up to the weight 
$ m = ( p^3 + p^2 + 2 ) -1 $; this was also our first example. In the 
lowest $ p^2 $-alcove the two algorithms agree, so we jump directly to the 
largest weight in that alcove.

\medskip
\medskip

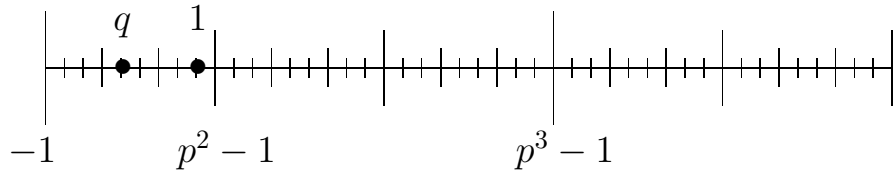
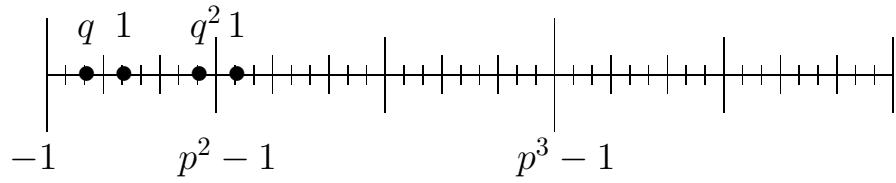
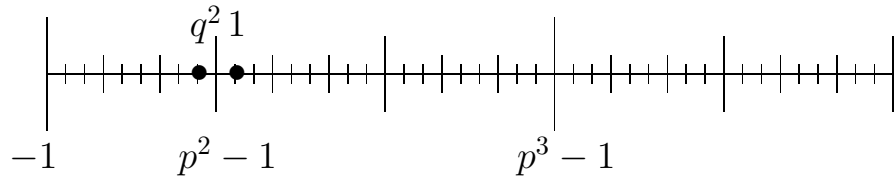
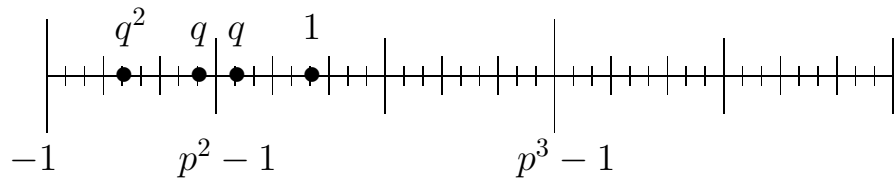
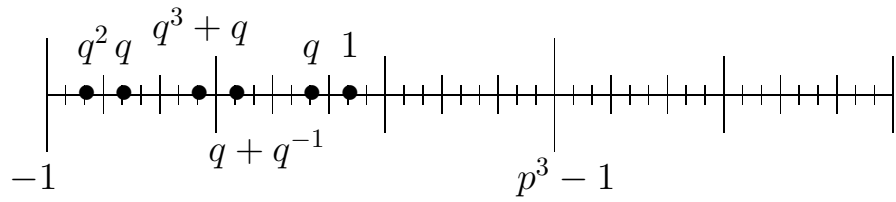
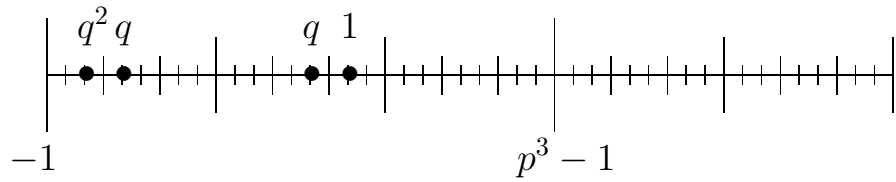
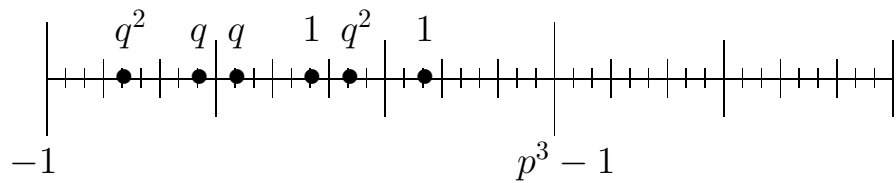
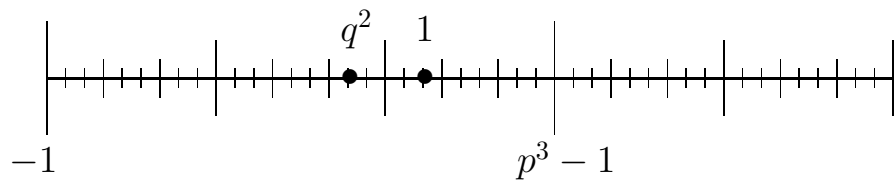
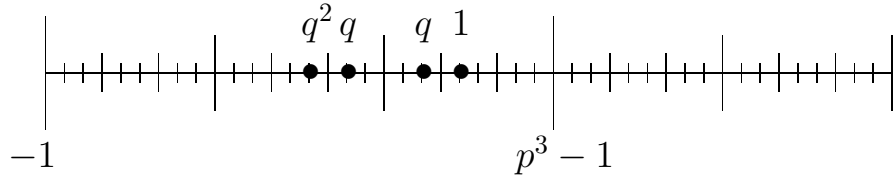
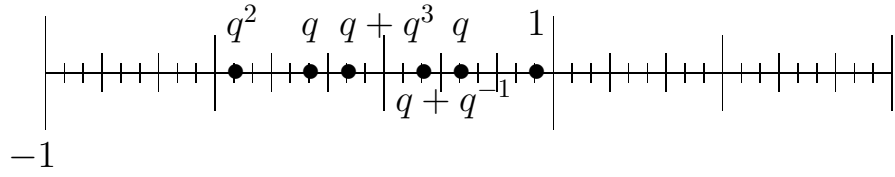
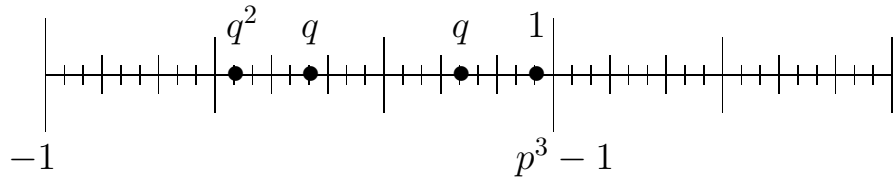
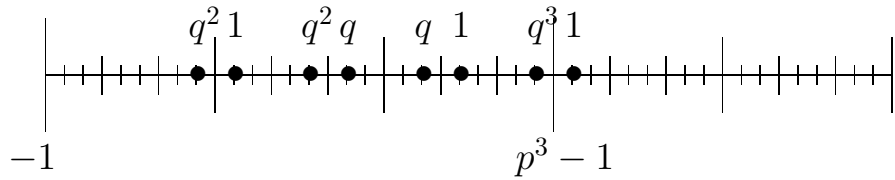
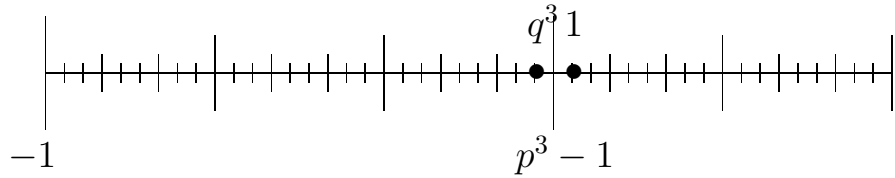
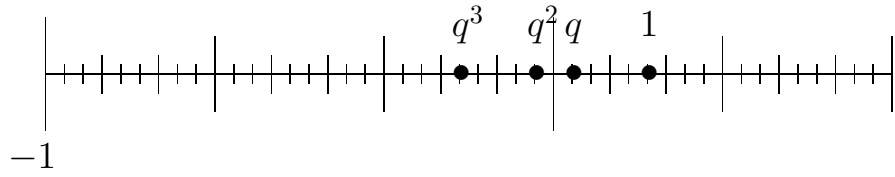
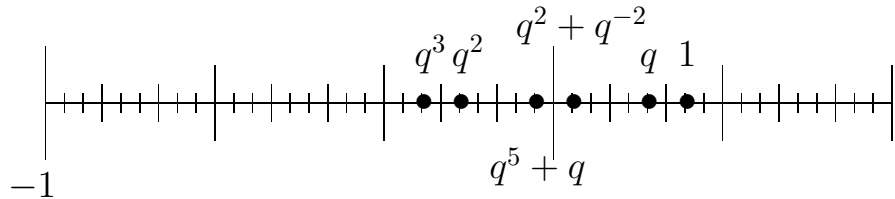
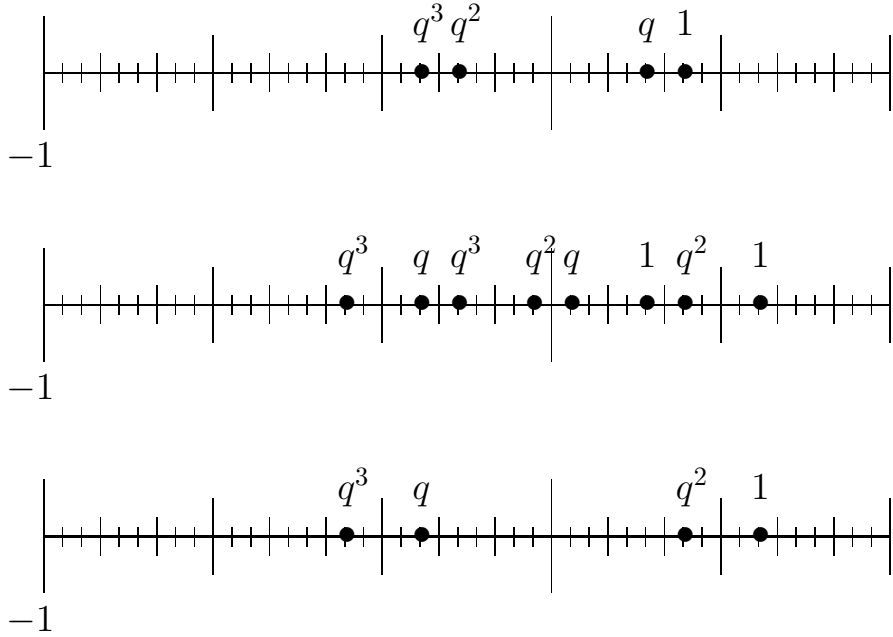
\setlength{\unitlength}{0.25cm}
\begin{center}
\begin{picture}(5,5)
\put(-20,0){\line(1,0){45}}
\linethickness{0.075mm}
  
\multiput(-20,-0.5)(1,0){46}{\line(0,1){1}}
  \multiput(-20,-1)(3,0){16}{\line(0,1){2}} 
  \multiput(-20,-2)(9,0){6}{\line(0,1){4}}
  \multiput(-20,-3)(27,0){2}{\line(0,1){6}} 
  \put(-22,-5){$-1$}
  \put(5,-5){$p^3-1$}
  \put(-13,-5){$p^2-1$}
  \put(-12.4,-0.4){$\bullet$}
  \put(-16.4,-0.4){$\bullet$}
 
  \put(-12.4,2){$1$}
  \put(-16.4,2){$q$}
 
\end{picture}
\end{center}

\medskip
\medskip
\medskip
\medskip
\medskip
\medskip
\medskip
\medskip

\setlength{\unitlength}{0.25cm}
\begin{center}
\begin{picture}(5,5)
\put(-20,0){\line(1,0){45}}
\linethickness{0.075mm}
  
\multiput(-20,-0.5)(1,0){46}{\line(0,1){1}}
  \multiput(-20,-1)(3,0){16}{\line(0,1){2}} 
  \multiput(-20,-2)(9,0){6}{\line(0,1){4}}
  \multiput(-20,-3)(27,0){2}{\line(0,1){6}} 
  \put(-22,-5){$-1$}
  \put(5,-5){$p^3-1$}
  \put(-13,-5){$p^2-1$}

  \put(-12.4,-0.4){$\bullet$}
  \put(-16.4,-0.4){$\bullet$}
  \put(-10.4,-0.4){$\bullet$}
  \put(-18.4,-0.4){$\bullet$}

  \put(-10.4,2){$1$}
  \put(-12.4,2){$q^2$}
  \put(-16.4,2){$1$}
  \put(-18.4,2){$q$}
\end{picture}
\end{center}

\medskip
\medskip
\medskip
\medskip
\medskip
\medskip
\medskip
\medskip

\setlength{\unitlength}{0.25cm}
\begin{center}
\begin{picture}(5,5)
\put(-20,0){\line(1,0){45}}
\linethickness{0.075mm}
  
\multiput(-20,-0.5)(1,0){46}{\line(0,1){1}}
  \multiput(-20,-1)(3,0){16}{\line(0,1){2}} 
  \multiput(-20,-2)(9,0){6}{\line(0,1){4}}
  \multiput(-20,-3)(27,0){2}{\line(0,1){6}} 
  \put(-22,-5){$-1$}
  \put(5,-5){$p^3-1$}
  \put(-13,-5){$p^2-1$}
  \put(-12.4,-0.4){$\bullet$}
  \put(-10.4,-0.4){$\bullet$}
  \put(-10.4,2){$1$}
  \put(-12.4,2){$q^2$}

\end{picture}
\end{center}

\medskip
\medskip
\medskip
\medskip
\medskip
\medskip
\medskip
\medskip

\setlength{\unitlength}{0.25cm}
\begin{center}
\begin{picture}(5,5)
\put(-20,0){\line(1,0){45}}
\linethickness{0.075mm}
  
\multiput(-20,-0.5)(1,0){46}{\line(0,1){1}}
  \multiput(-20,-1)(3,0){16}{\line(0,1){2}} 
  \multiput(-20,-2)(9,0){6}{\line(0,1){4}}
  \multiput(-20,-3)(27,0){2}{\line(0,1){6}} 
  \put(-22,-5){$-1$}
  \put(5,-5){$p^3-1$}
  \put(-13,-5){$p^2-1$}

  \put(-12.4,-0.4){$\bullet$}
  \put(-10.4,-0.4){$\bullet$}
  \put(-16.4,-0.4){$\bullet$}
  \put(-6.4,-0.4){$\bullet$}

  \put(-16.4,2){$q^2
$}
  \put(-12.4,2){$q$}
  \put(-10.4,2){$q$}
  \put(-6.4,2){$1$}
  
\end{picture}
\end{center}

\medskip
\medskip
\medskip
\medskip
\medskip
\medskip
\medskip

\setlength{\unitlength}{0.25cm}
\begin{center}
\begin{picture}(5,5)
\put(-20,0){\line(1,0){45}}
\linethickness{0.075mm}
  
\multiput(-20,-0.5)(1,0){46}{\line(0,1){1}}
  \multiput(-20,-1)(3,0){16}{\line(0,1){2}} 
  \multiput(-20,-2)(9,0){6}{\line(0,1){4}}
  \multiput(-20,-3)(27,0){2}{\line(0,1){6}} 
  \put(-22,-5){$-1$}
  \put(5,-5){$p^3-1$}

  \put(-12.4,-0.4){$\bullet$}
  \put(-10.4,-0.4){$\bullet$}
  \put(-16.4,-0.4){$\bullet$}
  \put(-6.4,-0.4){$\bullet$}
  \put(-4.4,-0.4){$\bullet$}
  \put(-18.4,-0.4){$\bullet$}

  \put(-18.4,2){$q^2$}
  \put(-16.4,2){$q$}
  \put(-14.4,3){$ q^3+ q $}
  \put(-11.4,-3.5){$ q + q^{-1}  $}
  \put(-6.4,2){$q$}
  \put(-4.4,2){$1$}

\end{picture}
\end{center}

\medskip
\medskip
\medskip
\medskip
\medskip
\medskip
\medskip

\setlength{\unitlength}{0.25cm}
\begin{center}
\begin{picture}(5,5)
\put(-20,0){\line(1,0){45}}
\linethickness{0.075mm}
  
\multiput(-20,-0.5)(1,0){46}{\line(0,1){1}}
  \multiput(-20,-1)(3,0){16}{\line(0,1){2}} 
  \multiput(-20,-2)(9,0){6}{\line(0,1){4}}
  \multiput(-20,-3)(27,0){2}{\line(0,1){6}} 
  \put(-22,-5){$-1$}
  \put(5,-5){$p^3-1$}

         \put(-16.4,-0.4){$\bullet$}
  \put(-6.4,-0.4){$\bullet$}
  \put(-4.4,-0.4){$\bullet$}
  \put(-18.4,-0.4){$\bullet$}

  \put(-18.4,2){$q^2$}
  \put(-16.4,2){$q$}
 
  \put(-6.4,2){$q$}
  \put(-4.4,2){$1$}

\end{picture}
\end{center}

\medskip
\medskip
\medskip
\medskip
\medskip
\medskip

\setlength{\unitlength}{0.25cm}
\begin{center}
\begin{picture}(5,5)
\put(-20,0){\line(1,0){45}}
\linethickness{0.075mm}
  
\multiput(-20,-0.5)(1,0){46}{\line(0,1){1}}
  \multiput(-20,-1)(3,0){16}{\line(0,1){2}} 
  \multiput(-20,-2)(9,0){6}{\line(0,1){4}}
  \multiput(-20,-3)(27,0){2}{\line(0,1){6}} 
  \put(-22,-5){$-1$}
  \put(5,-5){$p^3-1$}

  \put(-12.4,-0.4){$\bullet$}
  \put(-16.4,-0.4){$\bullet$}
  \put(-6.4,-0.4){$\bullet$}
  \put(-4.4,-0.4){$\bullet$}
  \put(-0.4,-0.4){$\bullet$}
  \put(-10.4,-0.4){$\bullet$}

  \put(-12.4,2){$q$}
  \put(-16.4,2){$q^2$}
  \put(-6.4,2){$1$}
  \put(-4.4,2){$q^2$}
  \put(-0.4,2){$1$}
  \put(-10.4,2){$q$}

\end{picture}
\end{center}

\medskip
\medskip
\medskip
\medskip
\medskip
\medskip
\medskip

\setlength{\unitlength}{0.25cm}
\begin{center}
\begin{picture}(5,5)
\put(-20,0){\line(1,0){45}}
\linethickness{0.075mm}
  
\multiput(-20,-0.5)(1,0){46}{\line(0,1){1}}
  \multiput(-20,-1)(3,0){16}{\line(0,1){2}} 
  \multiput(-20,-2)(9,0){6}{\line(0,1){4}}
  \multiput(-20,-3)(27,0){2}{\line(0,1){6}} 
  \put(-22,-5){$-1$}
  \put(5,-5){$p^3-1$}

  \put(-4.4,-0.4){$\bullet$}
  \put(-0.4,-0.4){$\bullet$}

  \put(-4.4,2){$q^2$}
  \put(-0.4,2){$1$}

\end{picture}
\end{center}

\medskip
\medskip
\medskip
\medskip
\medskip
\medskip
\medskip

\setlength{\unitlength}{0.25cm}
\begin{center}
\begin{picture}(5,5)
\put(-20,0){\line(1,0){45}}
\linethickness{0.075mm}
  
\multiput(-20,-0.5)(1,0){46}{\line(0,1){1}}
  \multiput(-20,-1)(3,0){16}{\line(0,1){2}} 
  \multiput(-20,-2)(9,0){6}{\line(0,1){4}}
  \multiput(-20,-3)(27,0){2}{\line(0,1){6}} 
  \put(-22,-5){$-1$}
  \put(5,-5){$p^3-1$}

  \put(-4.4,-0.4){$\bullet$}
  \put(-0.4,-0.4){$\bullet$}
  \put(-6.4,-0.4){$\bullet$}
  \put(1.6,-0.4){$\bullet$}

  \put(-4.4,2){$q$}
  \put(-0.4,2){$q$}
  \put(-6.4,2){$q^2$}
  \put(1.6,2){$1$}

\end{picture}
\end{center}

\medskip
\medskip
\medskip
\medskip
\medskip
\medskip
\medskip

\setlength{\unitlength}{0.25cm}
\begin{center}
\begin{picture}(5,5)
\put(-20,0){\line(1,0){45}}
\linethickness{0.075mm}
  
\multiput(-20,-0.5)(1,0){46}{\line(0,1){1}}
  \multiput(-20,-1)(3,0){16}{\line(0,1){2}} 
  \multiput(-20,-2)(9,0){6}{\line(0,1){4}}
  \multiput(-20,-3)(27,0){2}{\line(0,1){6}} 
  \put(-22,-5){$-1$}
  
  \put(-4.4,-0.4){$\bullet$}
  \put(-0.4,-0.4){$\bullet$}
  \put(-6.4,-0.4){$\bullet$}
  \put(1.6,-0.4){$\bullet$}
  \put(5.6,-0.4){$\bullet$}
  \put(-10.4,-0.4){$\bullet$}

  \put(-4.4,2){$q+q^3$}
  \put(-1.4,-2){$q+q^{-1}$}
  \put(-6.4,2){$q$}
  \put(1.6,2){$q$}
  \put(5.6,2){$1$}
  \put(-10.4,2){$q^2$}
\end{picture}
\end{center}

\medskip
\medskip
\medskip
\medskip
\medskip
\medskip
\medskip

\setlength{\unitlength}{0.25cm}
\begin{center}
\begin{picture}(5,5)
\put(-20,0){\line(1,0){45}}
\linethickness{0.075mm}
  
\multiput(-20,-0.5)(1,0){46}{\line(0,1){1}}
  \multiput(-20,-1)(3,0){16}{\line(0,1){2}} 
  \multiput(-20,-2)(9,0){6}{\line(0,1){4}}
  \multiput(-20,-3)(27,0){2}{\line(0,1){6}} 
  \put(-22,-5){$-1$}
  \put(5,-5){$p^3-1$}

  \put(-6.4,-0.4){$\bullet$}
  \put(1.6,-0.4){$\bullet$}
  \put(5.6,-0.4){$\bullet$}
  \put(-10.4,-0.4){$\bullet$}

  \put(-6.4,2){$q$}
  \put(1.6,2){$q$}
  \put(5.6,2){$1$}
  \put(-10.4,2){$q^2$}
\end{picture}
\end{center}

\medskip
\medskip
\medskip
\medskip
\medskip
\medskip
\medskip

\setlength{\unitlength}{0.25cm}
\begin{center}
\begin{picture}(5,5)
\put(-20,0){\line(1,0){45}}
\linethickness{0.075mm}
  
  \multiput(-20,-0.5)(1,0){46}{\line(0,1){1}}
  \multiput(-20,-1)(3,0){16}{\line(0,1){2}} 
  \multiput(-20,-2)(9,0){6}{\line(0,1){4}}
  \multiput(-20,-3)(27,0){2}{\line(0,1){6}} 
  \put(-22,-5){$-1$}
  \put(5,-5){$p^3-1$}

  \put(-6.4,-0.4){$\bullet$}
  \put(1.6,-0.4){$\bullet$}
  \put(5.6,-0.4){$\bullet$}
  \put(-10.4,-0.4){$\bullet$}
  \put(7.6,-0.4){$\bullet$}
  \put(-0.4,-0.4){$\bullet$}
  \put(-4.4,-0.4){$\bullet$}
  \put(-12.4,-0.4){$\bullet$}

  \put(-6.4,2){$q^2$}
  \put(1.6,2){$1$}
  \put(5.6,2){$q^3$}
  \put(-10.4,2){$1$}
  \put(7.6,2){$1$}
  \put(-0.4,2){$q$}
  \put(-4.4,2){$q$}
  \put(-12.4,2){$q^2$}

\end{picture}
\end{center}

\medskip
\medskip
\medskip
\medskip
\medskip
\medskip
\medskip

\setlength{\unitlength}{0.25cm}
\begin{center}
\begin{picture}(5,5)
\put(-20,0){\line(1,0){45}}
\linethickness{0.075mm}
  
  \multiput(-20,-0.5)(1,0){46}{\line(0,1){1}}
  \multiput(-20,-1)(3,0){16}{\line(0,1){2}} 
  \multiput(-20,-2)(9,0){6}{\line(0,1){4}}
  \multiput(-20,-3)(27,0){2}{\line(0,1){6}} 
  \put(-22,-5){$-1$}
  \put(5,-5){$p^3-1$}
 
  \put(5.6,-0.4){$\bullet$}
  \put(7.6,-0.4){$\bullet$}
   
  \put(5.6,2){$q^3$}
  \put(7.6,2){$1$}

\end{picture}
\end{center}

\medskip
\medskip
\medskip
\medskip
\medskip
\medskip
\medskip

\setlength{\unitlength}{0.25cm}
\begin{center}
\begin{picture}(5,5)
\put(-20,0){\line(1,0){45}}
\linethickness{0.075mm}
  
  \multiput(-20,-0.5)(1,0){46}{\line(0,1){1}}
  \multiput(-20,-1)(3,0){16}{\line(0,1){2}} 
  \multiput(-20,-2)(9,0){6}{\line(0,1){4}}
  \multiput(-20,-3)(27,0){2}{\line(0,1){6}} 

  \put(-22,-5){$-1$}

  \put(5.6,-0.4){$\bullet$}
  \put(7.6,-0.4){$\bullet$}
  \put(11.6,-0.4){$\bullet$}
  \put(1.6,-0.4){$\bullet$}

  \put(5.6,2){$q^2$}
  \put(7.6,2){$q$}
  \put(11.6,2){$1$}
  \put(1.6,2){$q^3$}

\end{picture}
\end{center}

\medskip
\medskip
\medskip
\medskip
\medskip
\medskip
\medskip

\setlength{\unitlength}{0.25cm}
\begin{center}
\begin{picture}(5,5)
\put(-20,0){\line(1,0){45}}
\linethickness{0.075mm}
  
  \multiput(-20,-0.5)(1,0){46}{\line(0,1){1}}
  \multiput(-20,-1)(3,0){16}{\line(0,1){2}} 
  \multiput(-20,-2)(9,0){6}{\line(0,1){4}}
  \multiput(-20,-3)(27,0){2}{\line(0,1){6}} 

  \put(-22,-5){$-1$}

  \put(5.6,-0.4){$\bullet$}
  \put(7.6,-0.4){$\bullet$}
  \put(11.6,-0.4){$\bullet$}
  \put(1.6,-0.4){$\bullet$}
  \put(13.6,-0.4){$\bullet$}
  \put(-0.4,-0.4){$\bullet$}

  \put(3.6,-4){$q^5+q$}
  \put(5,3.5){$q^2+q^{-2}$}
  \put(11.6,2){$q$}
  \put(1.6,2){$q^2$}
  \put(13.6,2){$1$}
  \put(-0.4,2){$q^3$}

\end{picture}
\end{center}

\medskip
\medskip
\medskip
\medskip
\medskip
\medskip
\medskip

\setlength{\unitlength}{0.25cm}
\begin{center}
\begin{picture}(5,5)
\put(-20,0){\line(1,0){45}}
\linethickness{0.075mm}
  
  \multiput(-20,-0.5)(1,0){46}{\line(0,1){1}}
  \multiput(-20,-1)(3,0){16}{\line(0,1){2}} 
  \multiput(-20,-2)(9,0){6}{\line(0,1){4}}
  \multiput(-20,-3)(27,0){2}{\line(0,1){6}} 

  \put(-22,-5){$-1$}
 
  \put(11.6,-0.4){$\bullet$}
  \put(1.6,-0.4){$\bullet$}
  \put(13.6,-0.4){$\bullet$}
  \put(-0.4,-0.4){$\bullet$}

  \put(11.6,2){$q$}
  \put(1.6,2){$q^2$}
  \put(13.6,2){$1$}
  \put(-0.4,2){$q^3$}

\end{picture}
\end{center}

\medskip
\medskip
\medskip
\medskip
\medskip
\medskip
\medskip

\setlength{\unitlength}{0.25cm}
\begin{center}
\begin{picture}(5,5)
\put(-20,0){\line(1,0){45}}
\linethickness{0.075mm}
  
  \multiput(-20,-0.5)(1,0){46}{\line(0,1){1}}
  \multiput(-20,-1)(3,0){16}{\line(0,1){2}} 
  \multiput(-20,-2)(9,0){6}{\line(0,1){4}}
  \multiput(-20,-3)(27,0){2}{\line(0,1){6}} 

  \put(-22,-5){$-1$}
 
  \put(11.6,-0.4){$\bullet$}
  \put(1.6,-0.4){$\bullet$}
  \put(13.6,-0.4){$\bullet$}
  \put(-0.4,-0.4){$\bullet$}
  \put(17.6,-0.4){$\bullet$}
  \put(7.6,-0.4){$\bullet$}
  \put(-4.4,-0.4){$\bullet$}
  \put(5.6,-0.4){$\bullet$}

  \put(11.6,2){$1$}
  \put(1.6,2){$q^3$}
  \put(13.6,2){$q^2$}
  \put(-0.4,2){$q$}
  \put(17.6,2){$1$}
  \put(7.6,2){$q$}
  \put(-4.4,2){$q^3$}
  \put(5.6,2){$q^2$}

\end{picture}
\end{center}

\medskip
\medskip
\medskip
\medskip
\medskip
\medskip
\medskip

\setlength{\unitlength}{0.25cm}
\begin{center}
\begin{picture}(5,5)
\put(-20,0){\line(1,0){45}}
\linethickness{0.075mm}
  
  \multiput(-20,-0.5)(1,0){46}{\line(0,1){1}}
  \multiput(-20,-1)(3,0){16}{\line(0,1){2}} 
  \multiput(-20,-2)(9,0){6}{\line(0,1){4}}
  \multiput(-20,-3)(27,0){2}{\line(0,1){6}} 

  \put(-22,-5){$-1$}
 
  \put(13.6,-0.4){$\bullet$}
  \put(-0.4,-0.4){$\bullet$}
  \put(17.6,-0.4){$\bullet$}
  \put(-4.4,-0.4){$\bullet$}
 
  \put(13.6,2){$q^2$}
  \put(17.6,2){$1$}
  \put(-4.4,2){$q^3$}
  \put(-0.4,2){$q$}

\end{picture}
\end{center}

\medskip
\medskip
\medskip
\medskip
\medskip
\medskip
\medskip
\medskip
\medskip
\medskip

So for all weights smaller than our chosen $ m $ we get results that 
are compatible with Erdmann's Theorem.

\medskip
{\bf Acknowledgement}
The results of section three were obtained 
independently by G. Murphy. 
I wish to thank him for a useful conversation 
on the matters of section four.

\end{document}